\documentclass[numbook, envcountsect, envcountsame, envcountreset,
 smallextended]{svjour3}                  %
\usepackage{amsmath}
\usepackage{graphicx}
\usepackage{mathptmx}
\usepackage{latexsym}
\usepackage{amsfonts}

\numberwithin{equation}{section}
\smartqed
\newcommand{\cB}{{\mathcal B}}

\newcommand{\cF}{{\mathcal F}}

\newcommand{\te}{{\theta}}

\newcommand{\Om}{{\Omega}}

\newcommand{\ve}{{\varepsilon}}
\newcommand{\del}{{\delta}}

\newcommand{\gam}{{\gamma}}
\newcommand{\Gam}{{\Gamma}}

\newcommand{\sig}{{\sigma}}
\newcommand{\al}{{\alpha}}
\newcommand{\be}{{\beta}}
\newcommand{\ka}{{\kappa}}

\newcommand{\La}{{\Lambda}}

\newcommand{\bbN}{{\mathbb N}}
\newcommand{\bbP}{{\mathbb P}}
\newcommand{\bbR}{{\mathbb R}}

\newcommand{\bbZ}{{\mathbb Z}}
\newcommand{\bbI}{{\mathbb I}}

\begin{document}
\title{From PET to SPLIT.}
\titlerunning{ From PET to SPLIT} 
 \author{ Yuri Kifer}
\institute{
Institute of Mathematics,
The Hebrew University,
Jerusalem 91904, Israel;
\email{ kifer@math.huji.ac.il}}

\date{Received: date / Accepted: date}
\maketitle
\begin{abstract}\noindent
The polynomial ergodic theorem (PET) which appeared in \cite{Be} and
attracted substantial attention in ergodic theory studies the limits
of expressions having the form $1/N\sum_{n=1}^NT^{q_1(n)}f_1\cdots
T^{q_\ell (n)}f_\ell$ where $T$ is a weakly mixing measure
preserving transformation, $f_i$'s are bounded measurable functions
and $q_i$'s are polynomials taking on integer values on the
integers. Motivated partially by this result we obtain a central
limit theorem for expressions of the form \newline
\centerline{$1/\sqrt{N}\sum_{n=1}^N
(X_1(q_1(n))X_2(q_2(n))\cdots
 X_\ell(q_\ell(n))-a_1a_2\cdots 
a_\ell)$}
\newline (sum-product limit theorem--SPLIT) where $X_i$'s are fast
$\alpha$-mixing bounded stationary processes, $a_j=EX_j(0)$ and
$q_i$'s are positive functions taking on integer values on integers
with some growth conditions which are satisfied, for instance, when
$q_i$'s are polynomials of growing degrees. This result can be
applied to the case when $X_i(n)=T^nf_i$ where $T$ is a mixing
subshift of finite type, a hyperbolic diffeomorphism or an expanding
transformation taken with a Gibbs invariant measure, as well, as to
the case when $X_i(n)=f_i(\xi_n)$ where $\xi_n$ is a Markov chain
satisfying the Doeblin condition considered as a stationary process
with respect to its invariant measure.

\end{abstract}
\markboth{Y.Kifer}{SPLIT}
\renewcommand{\theequation}{\arabic{section}.\arabic{equation}}
\pagenumbering{arabic}

\section{Introduction}\label{sec1}\setcounter{equation}{0}

The polynomial ergodic theorem (PET) appeared in \cite{Be} sais
 that in the $L^2$-sense\newline
$\lim_{N\to\infty}1/N\sum_{n=1}^NT^{q_1(n)}f_1\cdots T^{q_\ell (n)}f_\ell=
\prod_{i=1}^\ell\int f_id\mu$ where $T$ is a measure $\mu$
preserving weakly mixing transformation, $f_i$'s are bounded
measurable functions and $q_i$'s are polynomials taking on integer
values on the integers and satisfying $q_{i+1}(n)-q_{i}(n)\to\infty$
as $n\to\infty,\, i=1,...,\ell -1$. 
This and related results (see, for instance, \cite{FW}, \cite{FK} and
references there) where motivated originally by the study of multiple
recurrence for dynamical systems.
Namely, if $f_i=\bbI_{A_i},\, i=1,...,\ell$ are indicators of some
measurable sets $A_i$ of positive measure $\mu$ then PET implies
that for $\mu$-almost all (a.a.) $x$ the event $\cap_{i=1}^\ell \{
T^{q_i(n)}x\in A_i\}$ occurs with the frequency
$\prod_{i=1}^\ell\mu(A_i)$, in particular, infinitely often.

The probability theory name for the ergodic theorem is the law of
large numbers and after verifying it the next natural question to
ask is whether a central limit theorem type result holds also true
in this framework though, as usual, under somewhat stronger
assumptions. In this paper we will obtain convergence in
distribution to the normal law as $N\to\infty$ of expressions having
the form 
\[
\frac 1{\sqrt N}\sum_{n=1}^N \big (
X_1(q_1(n))X_2(q_2(n))\cdots X_\ell(q_\ell(n))-\prod_{i=1}^\ell
a_i\big )
\]
 (sum--product limit theorem: SPLIT) where $a_i=EX_i(0)$,
$X_i$'s are exponentially fast $\al$-mixing bounded stationary
processes and $q_i$'s are positive increasing for large $n$
functions taking on integer values on the integers with some growth
conditions which are satisfied, for instance, when $q_i$'s are
polynomials of increasing degrees. We observe that unlike PETs our
SPLITs do not require $q_i$'s to be polynomials, and so we obtain
also some new sum--product ergodic theorems paying the price of much
stronger mixing assumptions than in PETs. As in other cases with
central limit theorem our SPLIT describes, in particular,
fluctuations of the number of multiple recurrencies mentioned above
from its average frequency. In fact, we will derive a functional central 
limit theorem type extension of the above result.

Our results are applicable, for instance, to the case when
$X_i(n)=f_i(\xi_n)$ for bounded measurable $f_i$'s and a Markov
chain $\xi_n$ in a space $M$ satisfying the Doeblin condition (see
\cite{IL}) taken with its invariant measure $\mu$ which yields, in
particular, that for any measurable sets $A_i\subset M$ with
$\mu(A_i)>0,\, i=1,...,\ell$ if $N(n)$ is the number of events
$\cap_{i=1}^\ell\{\xi_{q_i(k)}\in A_i\}$ for $k$ running between 1
and $n$ then $n^{-1/2}(N(n)-\prod_{i=1}^\ell\mu(A_i))$ is
asymptotically normal. Our SPLITs seem to be new even when
$X_i(n),\, n=0,1,2,...$ are independent identically distributed
(i.i.d.) random variables though in this case the proof is much
easier and the result holds true in more general circumstances (see
Section \ref{sec5}). Another important class of processes satisfying
our conditions comes from dynamical systems where $X_i(n)=
f_i(T^nx)$ with $T$ being a topologically mixing subshift of finite
type or a $C^2$ expanding endomorphism or an Axiom A (in particular,
Anosov) (see \cite{Bo}) diffeomorphisms considered in a neighborhood
of an attractor taken with a Gibbs invariant measure. Some other
dynamical systems which fit our setup will be mentioned in the next
section. For a particular case of $Tx=\te x$ (mod 1), $\te>1,\,
x\in[0,1]$, polynomial $q_i$'s and fast approximable by
trigonometric polynomials $f_i$'s a corresponding central limit
theorem appears in \cite{Fu} whose specific setup allows application
of the Fourier analysis machinery.

Our methods are completely different from the ones in the ergodic
theory papers cited above and we rely on splitting the products into
weakly dependent factors (so SPLIT is not only an abbriviation here)
so that our main tool which is the inequality estimating the
difference between expectation of a product and a product of
expectations via the $\al$-mixing coefficient could be applied.
Observe that the martingale approximation methods which are popular
in modern proofs of the central limit theorem do not seem to work
(at least, directly) in our setup in view of strong dependencies
between past and future terms of sums here.

In writing of this paper I benefited from conversations with
V.Bergelson and B.Weiss who asked right questions and indicated to
me some references. Parts of the work were done during my visits to
the PennState and the Humboldt universities in Spring--Summer of
2008 in the framework of the Shapiro fellowship and the Humboldt
prize reinvitation programm, respectively, and I thank both
institutions for excellent working conditions and both foundations
for support.

\section{Preliminaries and main results}\label{sec2}\setcounter{equation}{0}

Our setup consists of $\ell$ bounded stationary processes
$X_1,X_2,...,X_\ell$, $|X_j(n)|\leq D<\infty,\, j=1,...,\ell;
n=0,1,...$ on a probability space $(\Om,\cF,P)$ and of a family of
$\sig$-algebras $\cF_{kl}\subset\cF,\, -\infty\leq k\leq
l\leq\infty$ such that $\cF_{kl}\subset\cF_{k'l'}$ if $k'\leq k$ and
$l'\geq l$. Given such family of $\sig$-algebras the $\al$-mixing
coefficient is defined by
\[
\al(n)=\sup_{k\geq
0}\sup_{A\in\cF_{-\infty,k},B\in\cF_{k+n,\infty}}|P(A\cap
B)-P(A)P(B)|,\, n\geq 0.
\]
Set also
\[
\beta_j(n)=\sup_{m\geq 0}E|X_j(m)-E(X_j(m)|\cF_{m-n,m+n})|.
\]
We assume that for some $\ka>0$,
\begin{equation}\label{2.1}
\al(n)+\max_{1\leq j\leq\ell}\be_j(n)\leq\ka^{-1}e^{-\ka n}.
\end{equation}
In what follows we can always consider $X(m)$ and $\cF_{kl}$ with
$m,k,l\geq 0$ only and just set formally in the above definitions
$\cF_{kl}=\cF_{kl}$ for $k<0$ and $l\geq 0$.

Next, let $q_1(n),q_2(n),...,q_\ell(n)$ be nonnegative functions
taking on integer values on the integers and such that $q_1(n)$ is
linear, i.e.,
\begin{equation}\label{2.2}
q_1(n)=rn+p\quad\mbox{for integer}\quad r>0,\, p\geq 0,
\end{equation}
and there exists $\gam\in(0,1)$ so that for all $n\geq n_0>1$,
\begin{equation}\label{2.3}
q_j(n+1)\geq q_j(n)+n^\gam,\,\, j=2,...,\ell
\end{equation}
and
\begin{equation}\label{2.4}
 q_{j+1}([n^{1-\gam}])\geq q_j(n)n^\gam,\,\, j=1,...,\ell-1.
\end{equation}
Observe that (\ref{2.3}) and (\ref{2.4}) are satisfied when $q_i$'s
are polynomials of positive degrees growing with $i$.

\begin{theorem}\label{thm2.1} Set $a_j=EX_j(0)$ and assume that
the above conditions (\ref{2.1})--(\ref{2.3}) on the processes $X_j$
and the functions $q_j$, $j=1,...,\ell$ hold true. Then, as $N\to\infty$,
\begin{equation}\label{2.5}
\frac 1{\sqrt{N}}\sum_{n=0}^N\big ( \prod_{j=1}^\ell
X_j(q_j(n))-\prod_{j=1}^\ell a_j\big ),
\end{equation}
 converges in distribution to a normal random variable with zero
mean and the variance
\begin{equation}\label{2.6}
\sig^2=\sig^2_\ell=EX_1^2(0)\big (\prod_{j=2}^\ell EX_j^2(0)-
\prod_{j=2}^\ell a_j^2\big )+\sig_1^2\prod_{j=2}^\ell a_j^2
\end{equation}
where 
\begin{eqnarray}\label{2.7}
&\sig^2_1=\lim_{N\to\infty}\frac 1NE\big (\sum_{n=1}^N(X_1(q_1(n))-a_1)
\big )^2\\
&=EX_1^2(0)-a_1^2+2\sum_{n=1}^\infty E\big((X_1(rn)-a_1)(X_1(0)-a_1)\big ),
\nonumber\end{eqnarray}
$\sig^2=\sig_1^2$ if $\ell=1$ and the last series in (\ref{2.7}) converges. 
Furthermore, $\sig=0$ if and only if either $X_j(0)=0$ almost surely (a.s.)
for some $j\geq 1$ or $X_j(0)=a_j$ a.s. for all $j\geq 2$ and $\sig_1=0$.
Finally, $\sig_1=0$ if and only if for all $m=0,1,2,...$,
\begin{equation}\label{2.7+}
X_1(rm+p)-a_1=U^{m+1}X-U^mX\quad\mbox{a.s.}
\end{equation}
where $U$ is the unitary operator associated with the stationary process
\newline
$\{ X_1(rm+p),\, m=0,1,2,...\}$ and $X$ belongs to the Hilbert space of
random variables with finite second moments which are measurable with
respect to the $\sig$-algebra generated by $\{ X_1(rm+p),\, m=0,1,2,...\}$
(see, for instance, \cite{IL}, Ch. 16).
\end{theorem}

Observe that since $EX_j^2(0)\geq a_j^2$ by the Cauchy--Schwarz inequality 
the last assertion of Theorem \ref{thm2.1} concerning $\sig=0$ follows from
(\ref{2.6}) and (\ref{2.7}) while the equivalence of $\sig_1=0$ and the
representation (\ref{2.7+}) is rather well known since it concerns the 
standard central limit theorem for 
\[
\frac 1{\sqrt N}\sum_{n=0}^N\big (X_1(rn+p)
-a_1\big ).
\]
 Still, for readers' convenience we recall the argument that
(\ref{2.7+}) follows from $\sig_1=0$ in Corollary \ref{cor3.7} while the
opposite implication is clear.

Note also that the case when $q_1(n)$ grows faster than linearly in
$n$ also fits our setup since we can take $X_1\equiv 1$ which would
mean that, in fact, we start with $X_2$ and $q_2$. In this case
\begin{equation}\label{2.8}
\sig^2=\prod_{j=1}^\ell EX_j^2(0)- \prod_{j=1}^\ell a_j^2
\end{equation}
and $\sig^2>0$ unless all $X_j$'s are constants with probability one.

In Section \ref{sec4+} we will extend Theorem \ref{thm2.1} to a more
general result where two linear functions $q_i$ are allowed. Namely,
set $q_0(n)=n$ and $q_j,\, j=1,2,...,\ell$ as above where $q_1$ is
given by (\ref{2.2}) with $r\geq 2$. We add another stationary process
$X_0$ with $X_0(n)\leq D$ for all $n$ and set $a_0=EX_0(0)$. Then we
have the following assertion.

\begin{theorem}\label{thm2.1+} As $N\to\infty$ the sequence of random
variables 
\begin{equation}\label{2.9+}
\frac 1{\sqrt{N}}\sum_{n=0}^N\big ( \prod_{j=0}^\ell
X_j(q_j(n))-\prod_{j=0}^\ell a_j\big ),
\end{equation}
converges in distribution to a normal random variable with zero mean and
the variance 
\begin{equation}\label{2.9++}
\hat\sig^2=EX_0^2(0)EX_1^2(0)\big (\prod_{j=2}^\ell EX_j^2(0)-
\prod_{j=2}^\ell a_j^2\big )+\sig_{01}^2\prod_{j=2}^\ell a_j^2
\end{equation}
where 
\begin{eqnarray}\label{2.9+++}
&\sig^2_{01}=\lim_{N\to\infty}\frac 1NE\big (\sum_{n=1}^N(X_0(n)X_1(q_1(n))
-a_0a_1)\big )^2\\
&=(EX_0^2(0)-a_0^2)EX_1^2(0)+a_0^2(EX_1^2(0)-a_1^2)
+2\sum_{n=1}^\infty E\big((X_0(n)-a_0)(X_0(0)-a_0)\big )\nonumber\\
&+2a_0^2\sum_{n=1}^\infty E\big((X_1(rn)-a_1)(X_1(0)-a_1)\big )+\Xi\nonumber
\nonumber\end{eqnarray}
and
\[
\Xi=2a_0a_1E\big (\sum_{n=0}^\infty (X_0(n)-a_0)(X_1(0)-a_1)
+\sum_{n=1}^\infty (X_0(0)-a_0)(X_1(n)-a_1)\big ).
\]
\end{theorem}

If we take $X_0\equiv 1$ then Theorem \ref{thm2.1+} reduces to Theorem
\ref{thm2.1+} where we need only $r\geq 1$. Furthermore, we can take
instead $X_j\equiv 1$ for all $j\geq 2$ which yields a nontrivial particular
case of Theorem \ref{thm2.1+} saying that 
\[
\frac 1{\sqrt{N}}\sum_{n=0}^N\big (X_0(n)X_1(rn+p)-a_0a_1\big ),\,\, r\geq 2,
\,\, r,p\in\bbN
\]
is asymptotically normal. 

For the readers' sake we will present first a complete proof of Theorem 
\ref{thm2.1} and then in Section \ref{sec4+} we explain additional elements
of the proof needed for Theorem \ref{thm2.1+} since a direct exposition from
the beginning of the latter more general case would make the reading more
difficult. Our main tool is splitting the products of $X_j(q_j(n_i))-\tilde
a_j$, where $\tilde a_j=0$ or $\tilde a_j=a_j$, in the way which
enables us to replace the expectation of a product by a product of
expectations with a sufficiently small error which will yield,
first, Gaussian type moment estimates for the expression in
(\ref{2.5}). Then we break the whole sum into a sum of blocks plus
terms which can be disregarded but play the role of gaps between
blocks. This will enable us to replace the characteristic function
of a sum of these blocks by a product of their characteristic
functions making only a small error. This is a standard method of
proving central limit theorem type results when such blocks can be
made sufficiently weakly dependent but in our case the terms of sums
depend on the far away future so our blocks are strongly dependent
and still, somewhat surprisingly, using the Taylor expansion of
characteristic functions and splitting products as described above
we can rely on this method in our case, as well. We observe that in
the case of Theorem \ref{thm2.1+} we will need, in fact, certain 
sequences of blocks so that the numbers $q_1(n),\, q_1(q_1(n)),\,
q_1(q_1(q_1(n))),...$ stay within the same sequence.  

Our $\al$-mixing condition is formulated in the form which allow
functions depending on the whole path of a stochastic process and
the exponentially fast decay (\ref{2.1}) holds true for many
important models. Let, for instance, $\xi_n$ be a Markov chain on a
space $M$ satisfying the Doeblin condition (see, for instance,
\cite{IL}, p.p. 367--368) and $f_j,\, j=1,...,\ell$ be a bounded measurable
functions
on the space of sequences $x=(x_i,\, i=0,1,2,...),\, x_i\in M$ such
that $|f_j(x)-f_j(y)|\leq Ce^{-cn}$ provided $x=(x_i),\, y=(y_i)$
and $x_i=y_i$ for all $i=0,1,...,n$ where $c,C>0$ do not depend on
$n$ and $j$. Set $X_j(n)=f_j(\xi_n,\xi_{n+1},\xi_{n+2},...)$ and let
$\sig$-algebras $\cF_{kl},\, k<l$ be generated by
$\xi_k,\xi_{k+1},...,\xi_l$ then the condition (\ref{2.1}) will be
satisfied considering $\{\xi_n,\, n\geq 0\}$ with its invariant
measure as a stationary process.

Important classes of processes satisfying our conditions come from
dynamical systems. Let $T$ be a $C^2$ Axiom A diffeomorphism (in
particular, Anosov) in a neighborhood of an attractor or let $T$ be
an expanding $C^2$ endomorphism of a Riemmanian manifold $M$ (see
\cite{Bo}), $f_j$'s are H\" older continuous functions and
$X_j(n)=f_j(T^nx)$. Here the probability space is $(M,\cB,\mu)$
where $\mu$ is a Gibbs invariant measure corresponding to some H\"
older continuous function. Let $\zeta$ be a finite Markov partition
for $T$ then we can take $\cF_{kl}$ to be the finite $\sig$-algebra
generated by the partition $\cap_{i=k}^lT^i\zeta$. In fact, we can
take here not only H\" older continuous $f_j$'s but also indicators
of sets from $\cF_{kl}$. A related example corresponds to $T$ being
a topologically mixing subshift of finite type which means that $T$
is the left shift on a subspace $\Xi$ of the space of one-sided
sequences $\xi=(\xi_i,i\geq 0), \xi_i=1,...,m$ such that $\xi\in\Xi$
if $\pi_{\xi_i\xi_{i+1}}=1$ where $\Pi=(\pi_{ij})$ is an $m\times m$
matrix with $0$ and $1$ entries and such that $\Pi^n$ for some $n$
is a matrix with positive entries. Again, we have to take in this
case $f_j$ to be H\" older continuous bounded functions of the
sequence space above, $\mu$ to be a Gibbs invariant measure
corresponding to some H\" older continuous function and to define
$\cF_{kl}$ as the finite $\sig$-algebra generated by cylinder sets
with fixed coordinates having numbers from $k$ to $l$. The
exponentially fast $\al$-(and even stronger)-mixing is well known in
the above cases (see \cite{Bo}). Among other dynamical systems with
exponentially fast $\al$-mixing we can mention also the Gauss map
$Tx=\{1/x\}$ of the unit interval with respect to the Gauss measure
(see \cite{He}).

A functional central limit theorem extension of Theorem \ref{thm2.1} can be 
derived by essentially the same method. Namely, for each $u\in[0,1]$ set
\begin{equation}\label{2.9}
W_N(u)=N^{-1/2}\sum_{n=0}^{[uN]}\big (\prod_{j=0}^\ell X_j(q_j(n))-
\prod_{j=0}^\ell a_j\big ).
\end{equation}
The process $W_N$ is a c\' adl\' ag, i.e. its paths belong to the space
$D[0,1]$ of right continuous functions on $[0,1]$ which have left limits
and, as usual, we consider $D[0,1]$ with the Skorokhod topology (see
\cite{Bi}). Denote by $W$ the standard one dimensional Brownian motion and
let $\bbP_{W_N}$ and $\bbP_{\hat\sig W}$ be the distributions of $W_N$ and 
of $\sig W(u),\, u\in[0,1]$ on $D[0,1]$, respectively, i.e.
\begin{equation}\label{2.10}
\bbP_{W_N}=P\{ W_N\in\Gam\}\quad\mbox{and}\quad \bbP_{\hat\sig W}(\Gam)=
P\{\hat\sig W\in\Gam\}
\end{equation}
for any Borel subset $\Gam$ of $D[0,1]$.

\begin{theorem}\label{thm2.2}
Under the conditions of Theorem \ref{thm2.1},
\begin{equation}\label{2.11}
\bbP_{W_N}\Rightarrow\bbP_{\hat\sig W}\quad\mbox{as}\quad N\to\infty
\end{equation}
where $\Rightarrow$ denotes the weak convergence of measures.
\end{theorem}

We will derive in Section \ref{sec4} Theorem \ref{thm2.2}, first, for the 
setup of Theorem \ref{thm2.1}, i.e. when $X_0\equiv 1$ and $\hat\sig=\sig$,
and the additional arguments of Section \ref{sec4+} will yield the result 
in the full generality of the setup of Theorem \ref{thm2.1+}. The proof
proceeds in the traditional way which consists of two ingredients.
First, we show by the block technique of Section \ref{sec4} (and by the
corresponding modification of Section \ref{sec4+}) that finite dimensional
distributions of $W_N$ weakly converge to corresponding finite dimensional 
distributions of $\hat\sig W$ which identifies the limit in (\ref{2.11}) 
uniquely (if it exists). Secondly, relying on Lemma \ref{lem3.7} (and its
generalisation in Section \ref{sec4+}) we obtain 
tightness of the family $\{\bbP_{W_N},\, N=1,2,...\}$ which yields the 
convergence.

\section{Gaussian type moment estimates}\label{sec3}\setcounter{equation}{0}

We start with the well known $\al$-mixing inequality (see, for
instance, \cite{Do} or \cite{DD}) saying that for any nonnegative
integers $k,n$ and random variables $Y$ and $Z$ which are
$\cF_{-\infty,k}$- and $\cF_{k+n,\infty}$-measurable, respectively,
\begin{equation}\label{3.1}
|E(YZ)-EYEZ|\leq 4\al(n)\| Y\|_\infty\| Z\|_\infty
\end{equation}
where $\|\cdot\|_\infty$ is the $L^\infty$-norm. This inequality
yields the following "splitting" lemma which will be our main
working tool throughout this paper.

\begin{lemma}\label{lem3.1}
Let $Y(j),\, j=0,1,...$ be bounded random variables and set
\begin{equation}\label{3.2}
\be(n)=\sup_{j\geq 0} E\big\vert
Y(j)-E(Y(j)|\cF_{j-n,j+n})\big\vert.
\end{equation}
Then for any $0\leq n_1\leq ...\leq n_l<n_{l+1}\leq n_{l+2}\leq
...\leq n_m$,
\begin{eqnarray}\label{3.3}
&\big\vert
E\prod_{i=1}^mY(n_i)-E\prod_{i=1}^lY(n_i)E\prod_{i=l+1}^mY(n_i)\big\vert\\
&\leq 2(m\be(k)+2\al(k))\prod_{i=1}^m\max(1,\| Y(n_i)\|_\infty)\nonumber
\end{eqnarray}
where $k=[(n_{l+1}-n_l)/3]$ and $[\cdot]$ denotes the integral part.
\end{lemma}
\begin{proof} Clearly,
\begin{eqnarray}\label{3.4}
&\big\vert
E\prod_{i=1}^mY(n_i)-E\prod_{i=1}^lY(n_i)E\prod_{i=l+1}^mY(n_i)\big\vert\\
&\leq I_1\prod_{i=l+1}^m\| Y(n_i)\|_\infty +I_2\prod_{i=1}^l\|
Y(n_i)\|_\infty +I_3\nonumber
\end{eqnarray}
where
\begin{equation}\label{3.5}
I_1=E\big\vert\prod_{i=1}^lY(n_i)-E(\prod_{i=1}^lY(n_i)|\cF_{-\infty,n_l+k})
\big\vert,
\end{equation}
\begin{equation}\label{3.6}
I_2=E\big\vert\prod_{i=l+1}^mY(n_i)-E(\prod_{i=l+1}^mY(n_i)|\cF_{n_{l+1}-k,
\infty})\big\vert
\end{equation}
and by (\ref{3.1}),
\begin{eqnarray}\label{3.7}
&I_3=\big\vert
E\big(E(\prod_{i=1}^lY(n_i)|\cF_{-\infty,n_l+k})E(\prod_{i=l+1}^mY(n_i)|
\cF_{n_{l+1}-k,\infty})\big
)\\
&-E\prod_{i=1}^lY(n_i)E\prod_{i=l+1}^mY(n_i)\big\vert\leq
4\al(k)\prod_{i=1}^m\| Y(n_i)\|_\infty.\nonumber
\end{eqnarray}
Observe that
\begin{eqnarray*}
&\big\vert
\prod_{i=1}^lY(n_i)-\prod_{i=1}^lE(Y(n_i)|\cF_{-\infty,n_l+k})\big\vert
\leq \sum_{j=1}^l\\
 &\big\vert\prod_{i=1}^{j-1}Y(n_i)\big
(Y(n_j)-E(Y(n_j)|\cF_{-\infty,n_l+k})\big
)\prod_{i=j+1}^lE(Y(n_i)|\cF_{-\infty,n_l+k})\big\vert\nonumber
\end{eqnarray*}
which together with (\ref{3.2}) and (\ref{3.5}) yields that
\begin{equation}\label{3.8}
I_1\leq 2l\be(k)\prod_{i=1}^l\max(1,\|Y(n_i)\|_\infty).
\end{equation}
Similarly,
\begin{equation}\label{3.9}
I_2\leq 2(m-l)\be(k)\prod_{i=l+1}^m\max(1,\|Y(n_i)\|_\infty),
\end{equation}
and so (\ref{3.3}) follows from (\ref{3.4})--(\ref{3.7}),
(\ref{3.8}) and (\ref{3.9}).
\qed\end{proof}

Next, set
\begin{eqnarray}\label{3.10}
&R(n)=\prod_{j=1}^\ell X_j(q_j(n))-\prod_{j=1}^\ell
a_j\\
&=\sum_{j=1}^\ell a_1\cdots
a_{j-1}(X_j(q_j(n))-a_j)X_{j+1}(q_{j+1}(n))\cdots X_\ell(q_\ell(n)).\nonumber
\end{eqnarray}
Here and in what follows if $\ell=1$ and a formula includes products of
undefined factors such as $X_{j+1},\, X_{j-1},\, a_{j+1},\, a_{j-1}$ with 
$1\leq j\leq\ell$ then such products should be replaced by $1$.
Observe that by (\ref{2.3}) and (\ref{2.4}) for any $j=1,...,\ell
-1$ and $n\geq n_0$,
\begin{equation}\label{3.11}
q_{i+1}(n)-q_i(n)\geq n-[n^{1-\gam}],
\end{equation}
and so by (\ref{3.3}) for such $n$,
\begin{equation}\label{3.12}
|ER(n)|\leq 4\ell D^\ell\big
(\ell\be([(n-[n^{1-\gam}])/3])+2\al([(n-[n^{1-\gam}])/3])\big ).
\end{equation}

The following result provides a Gaussian type estimate for the
second moment of sums of $R(n)$'s.

\begin{lemma}\label{lem3.2} There exists $C>0$ such that for all
$n\in\bbN$,
\begin{equation}\label{3.13}
E\big (\sum_{k=0}^nR(k)\big )^2\leq Cn.
\end{equation}
\end{lemma}
\begin{proof}
By (\ref{3.10}) for any $k_1,k_2\leq n$,
\begin{equation}\label{3.14}
|ER(k_1)R(k_2)|\leq\sum^\ell_{j_1,j_2=1}D^{j_1+j_2-2}|EQ_{j_1j_2}(k_1,k_2)|
\end{equation}
where, recall, $D$ is an upper bound on all $|X_j(k)|$'s and
\[
Q_{j_1j_2}(k_1,k_2)=\prod^2_{i=1}(X_{j_i}(q_{j_i}(k_i))-a_{j_i})X_{j_i+1}
(q_{j_i+1}(k_i))\cdots X_\ell(q_\ell(k_i)).
\]
Suppose that $q_{j_1}(k_1)<q_{j_2}(k_2)$ and $k_1,k_2>n_0$ where
$n_0$ is the same as in (\ref{2.3}) and (\ref{2.4}). Then by
(\ref{2.4}),
\[
q_{j_i}(k_i)<q_{j_i+1}(k_i)<...<q_\ell(k_i).
\]
Hence, we can apply (\ref{3.3}) with $k_i$ in place of $n_i$,
$Y(n_1)=X_{j_1}(q_{j_1}(k_1))-a_{j_1}$, $n_1=q_{j_1}(k_1)$, $l=1$
and $Y(n_i), i>1$ being other factors in the product for
$Q_{j_1j_2}(k_1,k_2)$ deriving that
\begin{equation}\label{3.15}
|EQ_{j_1j_2}(k_1,k_2)|\leq 16D^{2\ell}\big (
\ell\be(\nu_{j_1j_2}(k_1,k_2))+\al(\nu_{j_1j_2}(k_1,k_2))\big )
\end{equation}
where
\[
\nu_{j_1j_2}(k_1,k_2)=\min\big
([(q_{j_1+1}(k_1)-q_{j_1}(k_1))/3],[(q_{j_2}(k_2)-q_{j_1}(k_1))/3]\big
). \] This together with (\ref{2.1}), (\ref{2.3}) and (\ref{3.11})
yields that there exists a constant $C_1>0$ such that
\begin{equation}\label{3.16}
\sum_{1\leq j_1,j_2\leq\ell,n_0\leq k_1,k_2\leq n:\,
q_{j_1}(k_1)<q_{j_2}(k_2)}|EQ_{j_1j_2}(k_1,k_2)|\leq C_1(n+1).
\end{equation}
Now, if
\begin{equation}\label{3.17}
q_{j_1}(k_1)=q_{j_2}(k_2)
\end{equation}
for some $k_1,k_2\geq n_0$ then by (\ref{2.3}),
\[
q_{j_1}(k_1)<q_{j_2}(k_2+m)\quad\mbox{for all}\quad m\geq 1,
\]
and so the number of pairs $(j_2,k_2)$ such that $1\leq
j_2\leq\ell,\, n_0\leq k_2\leq n$ and (\ref{3.17}) is satisfied does
not exceed $\ell$. Hence, we obtain from here and (\ref{3.16}) that
\begin{equation}\label{3.18}
\sum_{1\leq j_1,j_2\leq\ell,n_0\leq k_1,k_2\leq
n}|EQ_{j_1j_2}(k_1,k_2)|\leq C_2(1+2\ell^2n_0+2\ell^2)(n+1)
\end{equation}
for some $C_2>0$ and (\ref{3.13}) follows.
\qed\end{proof}

\begin{remark}\label{rem3.3}
The estimates (\ref{3.11}) and (\ref{3.15}) enable us to obtain
(\ref{3.13}) under a weaker than (\ref{2.1}) condition, namely, a
polynomial decay of $\al(n)$ and $\be(n)$ so that either $
\sum_{n=1}^\infty(\al([n^\gam])+\be([n^\gam]))$ or
$\sum_{n=1}^\infty(\al([n^{1-\gam}])+\be([n^{1-\gam}]))$ converges
would already suffice. If we were interested only in (\ref{3.13}) we
could also weaken the boundedness condition on the stationary
processes $X_j,\, j=1,...,\ell$ assuming only existence of their
sufficiently high moments and using in place of (\ref{3.1}) the
inequality (see \cite{Do} or \cite{DD}),
\begin{equation}\label{3.19}
|E(YZ)-EYEZ|\leq 10\| Y\|_p\| Z\|_q(\al(n))^{1-\frac 1p-\frac 1q}
\end{equation}
which holds true provided $Y$ and $Z$ are $\cF_{-\infty,k}-$ and
$\cF_{k+n,\infty}-$measurable random variables, respectively, such
that $E|Y|^p<\infty,\, E|Z|^q<\infty$ and $\frac 1p+\frac 1q<1$.
Furthermore, (\ref{3.13}) does not require the full strength of the
assumption (\ref{2.4}) as we use only (\ref{3.11}) so that in place
of (\ref{2.4}) we can assume here, for instance, that $q_{j+1}(n)-
q_j(n)\geq\del n^\del$ for some $\del>0$ and all $n\geq n_0$.
\end{remark}
\begin{remark}\label{rem3.4}
Lemma \ref{lem3.2} yields that in the $L^2$-sense,
\begin{equation}\label{3.20}
\frac 1n\sum_{k=0}^nX_j(q_j(k))\longrightarrow\prod_{j=1}^\ell
a_j\,\,\mbox{as}\,\, n\to\infty
\end{equation}
which seems to be new when $q_j$'s are not polynomials.
\end{remark}

The following result justifies the formula (\ref{2.6}) for the
variance in our SPLIT.

\begin{lemma}\label{lem3.5}
Suppose that $N\geq n>m\geq [N^{1-\gam}]\geq n_0$. Then
\begin{equation}\label{3.21}
|E\big (\sum_{k=m+1}^nR(k)\big )^2-(n-m)\sig^2|\leq\hat C
\end{equation}
for some constant $\hat C>0$ independent of $n,m$ and $N$, where
$\sig$ is given by (\ref{2.6}).
\end{lemma}
\begin{proof} By (\ref{3.10})
\begin{eqnarray}\label{3.22}
& E\big (R(k_1)R(k_2)\big )=\sum^\ell_{j=1}a_1^2\cdots
a^2_{j-1}EQ_{jj}(k_1,k_2)\\
&+\sum_{\ell\geq j_2>j_1}a_1\cdots a_{j_1-1}a_1\cdots a_{j_2-1}\big
(EQ_{j_1j_2}(k_1,k_2)+EQ_{j_1j_2}(k_2,k_1)\big )\nonumber
\end{eqnarray}
where $Q_{j_1j_2}(k_1,k_2)$ is the same as in (\ref{3.14}). First,
we estimate $EQ_{jj}(k_1,k_2)$ for $j\geq 2$ and $k_1\ne k_2$, say,
when $k_2>k_1$. Assuming that $k_1\geq m$ it follows from
(\ref{2.3}) and (\ref{3.11}) that
\begin{equation}\label{3.23}
q_j(k_2)\geq q_j(k_1)+m^\gam\,\,\mbox{and}\,\, q_{j+1}(k_1)\geq
q_j(k_1)+m-[m^{1-\gam}],
\end{equation}
and so we can apply (\ref{3.15}) in order to obtain
\begin{equation}\label{3.24}
|EQ_{jj}(k_1,k_2)|\leq 16D^{2\ell}\big
(\ell\be(\rho_1(N))+\al(\rho_1(N))\big )
\end{equation}
where
\[
\rho_1(N)=\min\big ([[N^{1-\gam}]^\gam/3],\,
[([N^{1-\gam}]-[[N^{1-\gam}]^{1-\gam}])/3]\big )
\]
since $m\geq [N^{1-\gam}]$.

Next, if $j_2>j_1$ and $k,l\geq [N^{1-\gam}]$ then
\begin{equation}\label{3.25}
q_{j_2}(l)\geq q_{j_1}(k)N^\gam\geq q_{j_1}(k)+N^\gam -1,
\end{equation}
and so by (\ref{3.15}) we conclude that
\begin{equation}\label{3.26}
|EQ_{j_1j_2}(k,l)|\leq 16D^{2\ell}\big
(\ell\be(\rho_2(N))+\al(\rho_2(N))\big )
\end{equation}
where $\rho_2(N)=[(N^\gam -1)/3]$.

It remains to deal with the terms $Q_{jj}(k,k)$ and
$Q_{11}(k_1,k_2)$. Taking into account (\ref{3.25}) we apply
(\ref{3.3}) with $Y(n_1)=\big (X_j(q_j(k))-a_j\big )^2,\,
n_1=q_j(k),\, l=1$ and $Y(n_{j+i})=X^2_{j+i}(q_{j+i}(k)),\,
n_{j+i}=q_{j+i}(k),\, i=1,...,\ell -j$. It follows that
\begin{eqnarray}\label{3.27}
&|EQ_{jj}(k,k)-E\big (X_j(q_j(k))-a_j\big
)^2E\prod_{i=1}^{\ell-j}X^2_{j+i}(q_{j+i}(k))|\\
&\leq 16D^{2\ell}\big (\ell\be(\rho_2(N))+\al(\rho_2(N))\big
).\nonumber
\end{eqnarray}
Applying the same argument $\ell-j-1$ times to the expectation of
the product in (\ref{3.27}) and taking into account stationarity of
the processes $X_j$ we obtain that
\begin{eqnarray*}
&|EQ_{jj}(k,k)-E\big (X_j(0)-a_j\big
)^2\prod_{i=1}^{\ell-j}EX^2_{j+i}(0)|\\
&\leq 16(\ell+1)D^{2\ell}\big (\ell\be(\rho_2(N))+\al(\rho_2(N))\big
)
\end{eqnarray*}
and since $E(X_j(0)-a_j)^2=EX_j^2(0)-a_j^2$ it follows that
\begin{eqnarray}\label{3.28}
&|\sum_{j=1}^\ell a_1^2\cdots a_{j-1}^2EQ_{jj}(k,k)-\prod_{j=1}^\ell
EX^2_j(0)+\prod_{j=1}^\ell a_j^2|\\
&\leq 16\ell (\ell+1)D^{2\ell}\big
(\ell\be(\rho_2(N))+\al(\rho_2(N))\big ).\nonumber
\end{eqnarray}

Finally, in view of (\ref{2.2}) and (\ref{2.3}) for $k_2>k_1\geq 
[N^{1-\gam}]$ we obtain relying on (\ref{3.3}) similarly to the above that
\begin{eqnarray}\label{3.29}
&|EQ_{11}(k_1,k_2)-E\big (X_1(r(k_2-k_1))-a_1)(X_1(0)-a_1)\big )
\prod_{j=2}^{\ell}a^2_j|\\
&\leq 32D^{2\ell}\ell\big (\ell(\be(\rho_2(N))+\al(\rho_2(N))\big
).\nonumber
\end{eqnarray}
Again, by (\ref{3.3}) we have also
\begin{eqnarray*}
&|E\big (X_1(r(k_2-k_1))-a_1)(X_1(0)-a_1)\big )|\\
&\leq 16D^2\big (\be([r(k_2-k_1)/3]+ \al([ r(k_2-k_1)/3])\big ).
\end{eqnarray*}
This together with (\ref{2.1}) yields that for some constant $C_3>0$
independent of $n,m$ and $N$,
\begin{eqnarray}\label{3.30}
&|(n-m)\sum_{i=0}^\infty E\big (X_1(ri)-a_1)(X_1(0)-a_1)\big
)\\
&-\sum_{k_1=m+1}^n\sum_{i=0}^{n-k_1}E\big
(X_1(ri)-a_1)(X_1(0)-a_1)\big|\leq C_3.\nonumber
\end{eqnarray}
Collecting (\ref{3.24}), (\ref{3.26}) and (\ref{3.28})--(\ref{3.30})
we arrive at (\ref{3.21}) taking into account (\ref{2.1}) which
completes the proof of the lemma.
\qed\end{proof}

\begin{corollary}\label{cor3.6}
\begin{equation}\label{3.31}
\lim_{N\to\infty}\frac 1N E\big (\sum_{n=0}^NR(n)\big )^2=\sig^2
\end{equation}
and if $\sig=0$ then as $N\to\infty$ the expression (\ref{2.5}) 
converges to zero in distribution.
\end{corollary}
\begin{proof}
By (\ref{3.13}) for any $M<N$,
\begin{eqnarray}\label{3.32}
&|E\big (\sum_{n=0}^NR(n)\big )^2-E\big (\sum_{n=M}^NR(n)\big )^2|\\
&=|E\big (\sum_{n=0}^{M-1}R(n)\big )\big(\sum_{n=0}^NR(n)\nonumber\\
&+\sum_{n=M}^NR(n)\big )|\leq\sqrt 2\big (E\big
(\sum_{n=0}^NR(n)\big )^2\big )^{1/2} \big (E\big
(\sum_{n=0}^NR(n)\big )^2\nonumber\\
&+E\big (\sum_{n=M}^NR(n)\big )^2\big )^{1/2}\leq 2\sqrt 2\sqrt
{MN}\nonumber
\end{eqnarray}
and (\ref{3.31}) follows from (\ref{3.21}) and (\ref{3.22}) taking
$M=[N^{1-\gam}]+1$. If $\sig=0$ then (\ref{3.31}) together with the
Chebyshev inequality yields that as $N\to\infty$ the expression (\ref{2.5})
converges to zero in probability, and so in distribution, and in this case 
the main assertion of Theorem \ref{thm2.1} follows. 
\qed\end{proof}

\begin{corollary}\label{cor3.7}
If $\sig_1=0$ then
\begin{equation}\label{3.32+}
\sup_nE\big (\sum_{j=0}^n(X_1(rj+p))-a_1)\big )^2<\infty,
\end{equation}
and the representation (\ref{2.7+}) holds true.
\end{corollary}
\begin{proof} The inequality (\ref{3.32+}) follows from (\ref{2.7}), 
(\ref{3.29}) and (\ref{3.30}), and so by Theorem 18.2.2 from \cite{IL}
the representation (\ref{2.7+}) takes place. 
\qed\end{proof}

The following result gives the 4th moment Gaussian type estimate
needed to bound the error in the Taylor expansions of the
characteristic functions.

\begin{lemma}\label{lem3.7}
There exists $\tilde C>0$ such that whenever $N\geq n>m\geq
[N^{1-\gam}]\geq n_0$ then
\begin{equation}\label{3.33}
E\big (\sum_{k=m+1}^nR(k)\big )^4\leq\tilde C(n-m)^2.
\end{equation}
\end{lemma}
\begin{proof}
We have
\begin{equation}\label{3.34}
E\big (\sum_{k=m+1}^nR(k)\big
)^4\leq\sum_{k_1,k_2,k_3,k_4=m+1}^nA_{k_1k_2k_3k_4}
\end{equation}
where by (\ref{3.10}) for any $k_1,k_2,k_3,k_4$,
\begin{eqnarray}\label{3.35}
&A_{k_1k_2k_3k_4}=|E\big (R(k_1)R(k_2)R(k_3)R(k_4)\big
)|\\
&\leq\sum^\ell_{j_1,j_2,j_3,j_4=1}D^{j_1+j_2+j_3+j_4-4}|Q_{j_1j_2j_3j_4}
(k_1,k_2,k_3,k_4)|
\nonumber\end{eqnarray}
with
 \begin{eqnarray*}
&Q_{j_1j_2j_3j_4}(k_1,k_2,k_3,k_4)\\
&=E\prod_{i=1}^4\big ((
X_{j_i}(q_{j_i}(k_i))-a_{j_i})X_{j_i+1}(q_{j_i+1}(k_i))\cdots
X_{\ell}(q_{\ell}(k_i))\big ).
\end{eqnarray*}
In estimating the terms in the right hand side of (\ref{3.35}) we
assume without loss of generality that $j_1\leq j_2\leq j_3\leq
j_4$. If $j_1<j_2$ then taking into account that
$k_1,k_2,k_3,k_4>m\geq [N^{1-\gam}]$ we conclude relying on
(\ref{3.3}) and using (\ref{3.25}) similarly to (\ref{3.26}) that in
this case
\begin{equation}\label{3.36}
|Q_{j_1j_2j_3j_4}(k_1,k_2,k_3,k_4)|\leq 64D^{4\ell}\big
(2\ell\be(\rho_2(N))+\al(\rho_2(N))\big )
\end{equation}
with the same $\rho_2(N)$ as in (\ref{3.26}).

Next, consider the case $j_1=j_2<j_3=j_4$. Then
\begin{eqnarray}\label{3.37}
&Q_{j_1j_2j_3j_4}(k_1,k_2,k_3,k_4)\\
&=\prod_{i=1}^2\big
(X_{j_1}(q_{j_1}(k_i))-a_{j_1}\big )Z_1Z_2\prod_{i=3}^4\big
(X_{j_3}(q_{j_3}(k_i))-a_{j_3}\big )Z_3\nonumber
\end{eqnarray}
where $Z_1$ is the product of terms $X_j(q_j(k_i))$ with $i=1,2$ and
$j_1<j<j_3$, $Z_2$ is the product of terms $X_{j_3}(q_{j_3}(k_i))$
with $i=1,2$ and $Z_3$ is the product of terms $X_j(q_j(k_i))$ with
$j_3<j\leq\ell$ and $i=1,2,3,4$. Then employing 3 times (\ref{3.3})
and using again (\ref{3.25}) we obtain in this case that
\begin{eqnarray}\label{3.38}
&|Q_{j_1j_2j_3j_4}(k_1,k_2,k_3,k_4)-E\prod_{i=1}^2(X_{j_1}(q_{j_1}(k_i))-a_{j_1}\big
)EZ_1\\
&\times E\big(Z_2\prod_{i=3}^4(X_{j_3}(q_{j_3}(k_i))-a_{j_3})\big)EZ_3|
\nonumber \\
&\leq 192D^{4\ell}\big (2\ell\be(\rho_2(N))+\al(\rho_2(N))\big
)\nonumber
\end{eqnarray}
By (\ref{2.2}) and (\ref{2.3}) we see that for any $k,l\geq n_0$,
\begin{equation}\label{3.39}
|q_{j_i}(k)-q_{j_i}(l)|\geq |k-l|,
\end{equation}
and so we derive from (\ref{3.3}) that
\begin{equation}\label{3.40}
|E\prod_{i=1}^2(X_{j_1}(q_{j_1}(k_i))-a_{j_1})|\leq 16D^2\big
(\be([|k_1-k_2|/3])+\al([|k_1-k_2|/3])\big ).
\end{equation}
Applying the same argument twice we obtain also that
\begin{eqnarray}\label{3.41}
&|E\big (Z_2\prod_{i=3}^4\big (X_{j_3}(q_{j_3}(k_i))-a_{j_3}\big
)\big )-EZ_2E\prod_{i=3}^4\big (X_{j_3}(q_{j_3}(k_i))\\
&-a_{j_3}\big )|\leq 16D^{2\ell}\big
(\ell\be(\rho_3(k_1,k_2,k_3,k_4))+\al(\rho_3(k_1,k_2,k_3,k_4))\big
),\nonumber
\end{eqnarray}
where
\[
\rho_3(k_1,k_2,k_3,k_4)=\frac
13\min_{i_1=1,2;i_2=3,4}|k_{i_1}-k_{i_2}|
\]
and
\begin{equation}\label{3.42}
|E\prod_{i=3}^4\big (X_{j_3}(q_{j_3}(k_i))-a_{j_3}\big )|\leq
16D^2\big (\be([|k_3-k_4|/3])+\al([|k_3-k_4|/3])\big ).
\end{equation}

Next, if $j_1=j_2<j_3<j_4$ then we represent
$Q_{j_1j_2j_3j_4}(k_1,k_2,k_3,k_4)$ again in the form (\ref{3.37})
but now applying 3 times (\ref{3.3}) we obtain
\begin{eqnarray}\label{3.43}
&|Q_{j_1j_2j_3j_4}(k_1,k_2,k_3,k_4)-E\prod_{i=1}^2(X_{j_1}(q_{j_1}(k_i))-a_{j_1}\big
)EZ_1\\
&\times
E\big(Z_2(X_{j_3}(q_{j_3}(k_3))-a_{j_3})\big)E\big(Z_3(X_{j_4}(q_{j_4}(k_4))
-a_{j_4})\big)|\nonumber \\
& \leq 192D^{4\ell}\big (2\ell\be(\rho_2(N))+\al(\rho_2(N))\big
).\nonumber
\end{eqnarray}
Similarly to (\ref{3.41}) and (\ref{3.42}) it follows that
\begin{eqnarray}\label{3.44}
&|E\big (Z_2(X_{j_3}(q_{j_3}(k_3))-a_{j_3})\big )|\\
&\leq 8D^{2\ell-1}\big (\ell\be(\rho_4(k_1,k_2,k_3))+
\al(\rho_4(k_1,k_2,k_3))\big )\nonumber
\end{eqnarray}
where
\[
\rho_4(k_1,k_2,k_3)=\frac 13\min(|k_1-k_3|,|k_2-k_3|).
\]

Now, if $j_1=j_2=j_3<j_4$ then
\begin{equation}\label{3.45}
Q_{j_1j_2j_3j_4}(k_1,k_2,k_3,k_4)=\prod_{i=1}^3\big
(X_{j_1}(q_{j_1}(k_i))-a_{j_1}\big )Z_4
\end{equation}
where $Z_4$ is the product of $X_{j_4}(q_{j_4}(k_4))-a_{j_4}$ and
the terms of the form $X_j(q_j(k_i))$ with $\ell\geq j>j_1$ and
$i=1,2,3,4$. In this case by (\ref{3.3}) and (\ref{3.25}),
\begin{eqnarray}\label{3.46}
&|Q_{j_1j_2j_3j_4}(k_1,k_2,k_3,k_4)-E\prod_{i=1}^3\big
(X_{j_1}(q_{j_1}(k_i))-a_{j_1}\big )EZ_4|\\
&\leq 192D^{4\ell}\big (2\ell\be(\rho_2(N))+\al(\rho_2(N))\big
).\nonumber
\end{eqnarray}
Applying (\ref{3.3}) and (\ref{3.39}) we obtain that
\begin{eqnarray}\label{3.47}
&|E\prod_{i=1}^3\big (X_{j_1}(q_{j_1}(k_i))-a_{j_1}\big )|\\
&\leq 16D^3\big (3\be(\rho_5(k_1,k_2,k_3))+2\al(\rho_5(k_1,k_2,k_3))\big )
\nonumber\end{eqnarray}
where
\[
\rho_5(k_1,k_2,k_3))=\frac 16\big (\max
(k_1,k_2,k_3)-\min(k_1,k_2,k_3)\big ).
\]

Finally, in the case $j_1=j_2=j_3=j_4$ we can write
\begin{equation}\label{3.48}
Q_{j_1j_2j_3j_4}(k_1,k_2,k_3,k_4)=\prod_{i=1}^4\big
(X_{j_1}(q_{j_1}(k_i))-a_{j_1}\big )Z_5
\end{equation}
where $Z_5$ is the product of the terms $X_j(q_j(k_i))$ with
$\ell\geq j>j_1$ and $i=1,2,3,4$. Then by (\ref{3.3}) and
(\ref{3.25}) we have that
\begin{eqnarray}\label{3.49}
&|Q_{j_1j_2j_3j_4}(k_1,k_2,k_3,k_4)-E\prod_{i=1}^4\big
(X_{j_1}(q_{j_1}(k_i))-a_{j_1}\big )EZ_5|\\
&\leq 192D^{4\ell}\big (2\ell\be(\rho_2(N))+\al(\rho_2(N))\big ).
\nonumber
\end{eqnarray}
Suppose that $k_{i_1}\leq k_{i_2}\leq k_{i_3}\leq k_{i_4}$ where
$i_1,i_2,i_3,i_4$ are different integers between 1 and 4. Then by
(\ref{3.3}) and (\ref{3.25}),
\begin{eqnarray}\label{3.50}
&|E\prod_{l=1}^4\big(X_{j_1}(q_{j_1}(k_{i_l}))-a_{j_1}\big )|\\
&\leq 64D^4\big (2\be(\rho_6(k_1,k_2,k_3,k_4))+\al(\rho_6(k_1,k_2,k_3,k_4))
\big )\nonumber
\end{eqnarray}
where
\[
\rho_6(k_1,k_2,k_3,k_4)=\frac 13\max(|k_{i_2}-k_{i_1}|,
|k_{i_4}-k_{i_3}|).
\]
Collecting (\ref{3.34})--(\ref{3.38}) and (\ref{3.40})--(\ref{3.50})
and taking into account (\ref{2.1}) we arrive at (\ref{3.33})
completing the proof of the lemma.
\qed\end{proof}

\begin{remark}\label{rem3.8}
It is clear from the above arguments that the proofs of Lemmas
\ref{lem3.5} and \ref{lem3.7} still go through if in place of
(\ref{3.1}) and boundedness of $X_j$'s we assume that $\al(n)$ and
$\be(n)$ decay with sufficiently fast polynomial speed and some high
enough moments of $X_j$'s are finite so that we could apply
(\ref{3.19}) sufficiently many times. This would not suffice in the
next section where we have to apply (\ref{3.1}) in the form of
(\ref{3.3}) the number of times growing in $N$, and so (\ref{3.19})
with any fixed $p$ and $q$ will not work.
 \end{remark}

 \begin{remark}\label{rem3.9}
 Lemma \ref{lem3.7} yields the convergence (\ref{3.20}) with
 probability one. Indeed, (\ref{3.33}) together with Chebyshev's
 inequality gives that
 \[
 P\{\frac 1n|\sum_{k=0}^nR(k)|\geq\frac 1{n^8}\}\leq\tilde Cn^{-3/2}
 \]
 which in view of the Borel--Cantelly lemma implies the above
 assertion.
 \end{remark}

\section{Blocks and characteristic functions}\label{sec4}\setcounter{equation}{0}

Choose a small positive $\ve$ and a large $L\geq 4$ so that
$L\ve<\gam /4$. Set $\tau(N)=[N^{1-\ve}],\,\te(N)=[N^{1-L\ve}],\,
m(N)=\big [\frac N{\te(N)+\tau(N)}\big ]$ and introduce the sets of
integers
\[
\Gam_k(N)=\{ n:\,\te(N)+(k-1)(\te(N)+\tau(N))\leq n\leq
k(\te(N)+\tau(N))\}
\]
and
\[
\tilde\Gam_k(N)=\{ n:\, (k-1)(\te(N)+\tau(N))+1\leq n\leq
\te(N)+(k-1)(\te(N)+\tau(N))\}.
\]
Assuming that $N\geq\exp(2/\ve)$ which ensures that $m(N)\geq 1$ set
for $k=1,2,...,m(N)$,
\[
Y_k=\sum_{n\in\Gam_k(N)}R(n)\,\,\,\mbox{and}\,\,\,
Z_k=\sum_{n\in\tilde\Gam_k(N)}R(n)
\]
where $R(n)$ is the same as in (\ref{3.10}). Till the end of this
section our goal will be to show that the characteristic function
$\Phi_N(t)=E\exp\big (\frac {it}{\sqrt N}\sum_{n=0}^NR(n)\big )$
converges to $\exp(-\sig^2t^2/2)$ which will yield Theorem
\ref{thm2.1}. In doing so we employ the blocks (partial sums)
introduced above and the estimates of Section \ref{sec3} so that we
will deal mainly with the larger blocks $Y_k$ showing that the
smaller blocks $Z_k$ can be disregarded and they will be treated as
gaps between $Y_k$'s.

First, setting
\[
\Psi_N(t)=E\exp\big (\frac {it}{\sqrt N}\sum_{1\leq n\leq
m(N)}Y_n\big ) \]
 and relying on the inequality
\[
|e^{i(x+y)}-e^{iy}|=|e^{ix}-1|\leq |x|
\]
we obtain from (\ref{3.13}) and (\ref{3.33}) that
\begin{eqnarray}\label{4.1}
&|\Phi_N(t)-\Psi_N(t)|\leq\frac {|t|}{\sqrt N}E\big
(|\sum_{n=0}^{\te(N)}R(n)|+|\sum_{2\leq
n\leq m(N)}Z_n|\\
&+|\sum^N_{n=m(N)(\te(N)+\tau(N))+1}R(n)|\big )\leq\frac {|t|}{\sqrt
N}\big (\big (E(\sum_{n=0}^{\te(N)}R(n))^2\big )^{1/2}\nonumber\\
&+\sum_{2\leq n\leq m(N)}(EZ_n^4)^{1/4}+\big
(E(\sum^N_{n=m(N)(\te(N)+\tau(N))+1}R(n))^4\big )^{1/4}\big
)\nonumber\\
&\leq\frac {|t|}{\sqrt N}(\sqrt C\sqrt {\te(N)+1}+ \tilde
C^{1/4}m(N)\sqrt {\te(N)}+\tilde C^{1/4}\sqrt
{\te(N)+\tau(N)})\nonumber\\
&\leq\check C|t|(N^{-\ve(\frac L2-1)}+N^{-\ve/2})\nonumber
\end{eqnarray}
for some constant $\check C>0$ independent of $N$.

The main part of this section is the following result showing that
up to a small error the characteristic function of the sum of blocks
$Y_k$ is close to the product of characteristic functions of $Y_k$'s
themselves. When blocks are weakly dependent this step follows
immediately from (\ref{3.1}) but our blocks are strongly dependent,
and so the proof requires some work. Set
\[
\psi_N^{(k)}(t)=E\exp\big (\frac {it}{\sqrt N}Y_k\big ),\,\, k\leq
m(N).
\]

\begin{lemma}\label{lem4.1}
For any $t$ and each small $\ve>0$ there exists $K_\ve(t)>0$ such
that for all $N\geq\exp(2/\ve)$,
\begin{equation}\label{4.2}
|\Psi_N(t)-\prod_{1\leq k\leq m(N)}\psi_N^{(k)}(t)|\leq
K_\ve(t)N^{-\frac \ve 2\sqrt N}.
\end{equation}
\end{lemma}
\begin{proof} Set $\hat Y_k=Y_k+\tau(N)\prod_{j=1}^\ell a_j$,
\[
\hat\Psi_N(t)=E\exp\big (\frac {it}{\sqrt N}\sum_{1\leq k\leq
m(N)}\hat Y_k\big )\,\,\mbox{and}\,\,\hat\psi_N^{(k)}(t)=E\exp\big
(\frac {it}{\sqrt N}\hat Y_k\big ).
\]
Then, clearly,
\begin{equation}\label{4.3}
|\Psi_N(t)-\prod_{1\leq k\leq
m(N)}\psi_N^{(k)}(t)|=|\hat\Psi_N(t)-\prod_{1\leq k\leq
m(N)}\hat\psi_N^{(k)}(t)|.
\end{equation}
By the reminder formula for the Taylor expansion
\begin{equation}\label{4.4}
|e^{iz}-\sum^n_{k=0}\frac {(iz)^k}{k!}|\leq\frac
{|z|^{n+1}}{(n+1)!}.
\end{equation}
With the same $\ve>0$ as above set
\begin{equation}\label{4.5}
n(N)=n_\ve(N)=[N^{\frac 12+\ve}]
\end{equation}
and denote
\[
I_N^{(k)}(t)=\sum_{l=0}^{n(N)}\frac {(it)^l}{N^{l/2}l!}\hat Y_k^l.
\]
Then by (\ref{4.4}),
\begin{equation}\label{4.6}
|\exp\big (\frac {it}{\sqrt N}\hat Y_k\big )-I_N^{(k)}(t)|\leq\frac
{(|t|D\sqrt N)^{n(N)+1}}{(n(N)+1)!}\leq C_4^{n(N)}|t|^{n(N)}N^{-\ve
n(N)}
\end{equation}
for some constant $C_4>0$ independent of $N\geq 4$. Then
\begin{equation}\label{4.7}
|\hat\Psi_N(t)-\prod_{1\leq k\leq m(N)}\hat\psi_N^{(k)}(t)|\leq
J(t,N)+\del(t,N)
\end{equation}
where
\[
J(t,N)=|E\prod_{1\leq k\leq m(N)}I_N^{(k)}(t)-\prod_{1\leq k\leq
m(N)}EI_N^{(k)}(t)|
\]
and
\begin{eqnarray}\label{4.8}
&\del(t,N)=2m(N)C_4^{n(N)}|t|^{n(N)}N^{-\ve
n(N)}\\
&\times (1+C_4^{n(N)}|t|^{n(N)}N^{-\ve n(N)})^{m(N)}
\leq C(\ve,t)N^{-\frac \ve 2\sqrt N}\nonumber
\end{eqnarray}
for some $C(\ve,t)>0$ independent of $N$.

It remains to estimate $J(t,N)$ which is the main point of the
proof. We have
\begin{equation}\label{4.9}
J(t,N)=\sum_{0\leq l_1,...,l_{n(N)}\leq
n(N)}|tN^{-1/2}|^{\sum_{1\leq k\leq
m(N)}l_k}\prod_{k=1}^{m(N)}(l_k!)^{-1}G_{l_1,...,l_{m(N)}}(t,N)
\end{equation}
where
\[
G_{l_1,...,l_{m(N)}}(t,N)=|E\prod_{k=1}^{m(N)}\hat
Y_k^{l_k}-\prod_{k=1}^{m(N)}E\hat Y_k^{l_k}|.
\]
Next, we represent the $l_k$-th power of the sum $\hat Y_k$ in the
form
\begin{equation}\label{4.10}
Y_k^{l_k}=\sum_{\sig^{(k)}}\be_{\sig^{(k)}}^{(k)}\prod_{n\in\Gam_k(N)}
\prod_{j=1}^\ell X_j^{\sig_n^{(k)}}(q_j(n))
\end{equation}
where $\be^{(k)}_{\sig^{(k)}}$ are $l_k$-nomial coefficients and $
\sig^{(k)}=(\sig^{(k)}_n,\, n\in\Gam_k(N))$ satisfies
\begin{equation}\label{4.11}
\sig^{(k)}_n\geq
0\,\,\mbox{and}\,\,\sum_{n\in\Gam_k(N)}\sig^{(k)}_n=l_k\leq n(N).
\end{equation}
Then
\begin{equation}\label{4.12}
G_{l_1,...,l_{m(N)}}(t,N)\leq\sum_{\sig^{(1)},\sig^{(2)},...,\sig^{(m(N))}}
\prod_{k=1}^{m(N)}\be^{(k)}_{\sig^{(k)}}H_{l_1,...,l_{m(N)}}(t,N)
\end{equation}
where
\begin{eqnarray*}
&H_{l_1,...,l_{m(N)}}(t,N)=|E\prod_{k=1}^{m(N)}\prod_{n\in\Gam_k(N)}
\prod_{j=1}^\ell X_j^{\sig_n^{(k)}}(q_j(n))\\
&-\prod_{k=1}^{m(N)}E\prod_{n\in\Gam_k(N)}\prod_{j=1}^\ell
X_j^{\sig_n^{(k)}}(q_j(n))|.
\end{eqnarray*}

Next, we change the order of products in the two expectations above
so that the product $\prod_{j=1}^\ell$ appear immediately after the
expectation and apply the "splitting" Lemma \ref{lem3.1} $\ell$
times to the latter product for both expectations. Since $n\geq
[N^{1-L\ve}]$ in the above expressions then relying $\ell$ times on
(\ref{3.3}) and the second part of (\ref{3.23}) we obtain taking
into account (\ref{4.11}) that
\begin{eqnarray}\label{4.13}
&|E\prod_{k=1}^{m(N)}\prod_{n\in\Gam_k(N)}\prod_{j=1}^\ell
X_j^{\sig_n^{(k)}}(q_j(n))\\
& -\prod_{j=1}^\ell E\prod_{k=1}^{m(N)}\prod_{n\in\Gam_k(N)}
X_j^{\sig_n^{(k)}}(q_j(n))|\nonumber\\
&\leq 2\ell D^{\ell n(N)m(N)}\big (\ell
n(N)m(N)\be(\rho_6(N))+2\al(\rho_6(N))\big )\nonumber
\end{eqnarray}
where
\[
\rho_6(N)=[\frac 13([N^{1-L\ve}]-[N^{(1-\gam)(1-L\ve)}])].
\]
Similarly,
\begin{eqnarray}\label{4.14}
&|E\prod_{n\in\Gam_k(N)}\prod_{j=1}^\ell X_j^{\sig_n^{(k)}}(q_j(n))
-\prod_{j=1}^\ell E\prod_{n\in\Gam_k(N)}
X_j^{\sig_n^{(k)}}(q_j(n))|\\
&\leq 2\ell D^{\ell n(N)}\big (\ell
n(N)\be(\rho_6(N))+2\al(\rho_6(N))\big ).\nonumber
\end{eqnarray}

Next, for each fixed $j$ we apply (\ref{3.3}) $m(N)$ times to the
product $\prod_{k=1}^{m(N)}$ appearing after the expectation and in
view of (\ref{3.39}) and the size of the gaps $Z_k$ between the
blocks $Y_k$ it follows that
\begin{eqnarray}\label{4.15}
\,\,\,\,\,\, &|E\prod_{k=1}^{m(N)}\prod_{n\in\Gam_k(N)}
X_j^{\sig_n^{(k)}}(q_j(n)) -\prod_{k=1}^{m(N)}
E\prod_{n\in\Gam_k(N)}
X_j^{\sig_n^{(k)}}(q_j(n))|\\
&\leq 2m(N)D^{m(N)n(N)}\big (m(N)
n(N)\be([[N^{1-L\ve}]/3])+2\al([[N^{1-L\ve}]/3])\big ).\nonumber
\end{eqnarray}
Collecting (\ref{4.3}), (\ref{4.5})--(\ref{4.15}) and taking into
account that for each $k$,
\[
\sum_{\sig^{(k)}}\be^{(k)}_{\sig^{(k)}}\leq N^{(1-\ve)l_k}
\]
and
\[
\sum_{1\leq l_1,...,l_{m(N)}\leq n(N)}\prod_{k=1}^{m(N)}\frac
{|N^{\frac 12-\ve}t|^{l_k}}{l_k!}\leq\exp(N^{\frac 12-\ve}|t|m(N))
\]
we arrive at (\ref{4.2}).
\qed\end{proof}

Now we can complete the proof of Theorem \ref{thm2.1}. Using the
inequalities
\[
|e^{ix}-1-ix+\frac {x^2}{2}|\leq |x|^3\,\,\mbox{and}\,\,
|e^{-x}-1+x|\leq x^2
\]
which hold true for any real $x$ we derive from (\ref{3.12}),
(\ref{3.21}) and (\ref{3.33}) together with the H\" older inequality
that
\begin{eqnarray}\label{4.16}
&|\psi_N^{(k)}(t)-\exp\big (-\frac {\sig^2t^2\tau(N)}{2N}\big
)|\\
&\leq 4\ell D^\ell N^{\frac 12-\ve}|t|\big
(\ell\be(\rho_6(N))+2\al(\rho_6(N))\big )\nonumber \\
&+\tilde
C^{3/4}|t|^3N^{-3\ve/2}+\frac {\sig^4t^4}{4N^2}(\tau(N))^2\nonumber
\end{eqnarray}
where $\rho_6$ is the same as in (\ref{4.13}). Taking into account
that
\begin{equation}\label{4.17}
|\prod_{1\leq k\leq l}g_k-\prod_{1\leq k\leq l}h_k|\leq\sum_{1\leq
k\leq l}|g_k-h_k|
\end{equation}
whenever $0\leq |g_k|,|h_k|\leq 1,\, k=1,...,l$ we obtain from
(\ref{4.16}) that
 \begin{eqnarray}\label{4.18}
&|\prod_{1\leq k\leq m(N)}\psi_N^{(k)}(t)-\exp\big (-\frac
{\sig^2t^2}{2}\big)|\leq\frac {\sig^2t^2}{2}(1-\frac
{\tau(N)m(N)}{N})\\
&+4\ell D^\ell N^{\frac 12-\ve}m(N)|t|^3\big
(\ell\be(\rho_6(N))+2\al(\rho_6(N))\big )\nonumber \\
&+\tilde C^{3/4}|t|^3N^{-3\ve/2}m(N)+\frac
{\sig^4t^4}{4N^2}(\tau(N))^2m(N)\nonumber
\end{eqnarray}
and since $m(N)$ is of order $N^\ve$ while $\tau(N)$ is of order
$N^{1-\ve}$ we obtain that the right hand side of (\ref{4.18}) is
bounded by const$(t^4+1)N^{-\ve/2}$. This together with (\ref{4.1})
and (\ref{4.2}) gives
\begin{equation}\label{4.19}
|\Phi_N(t)-\exp(-\frac 12\sig^2t^2)|\leq\tilde K_\ve(t)N^{-\ve/2}
\end{equation}
for some $\tilde K_\ve(t)>0$ independent of $N$ and the assertion of
Theorem \ref{thm2.1} follows. \qed

Next, we explain the proof of Theorem \ref{thm2.2}. In order to show 
that finite dimensional distributions of $W_N$ converge to corresponding
finite dimensional distributions of $\sig W$ we fix $0=u_0<u_1<u_2<\cdots 
<u_k\leq 1$ and some real $t_1,t_2,...,t_k$ proving that 
\begin{eqnarray}\label{4.20}
&\Phi_N^{u_1,...,u_k}(t_1,...,t_k)\\
&=E\exp\big (i\sum_{j=1}^kt_jW_N(u_j)
\big )\longrightarrow\phi_{\sig W}^{u_1,...,u_k}(t_1,...,t_k)\nonumber\\
&=\prod_{j=1}^k\exp\big(-\frac 12\sig^2(u_j-u_{j-1})(\sum_{l=j}^kt_l)^2\big )
\,\,\,\mbox{as}\,\,\, N\to\infty .\nonumber
\end{eqnarray}
First, we have
\begin{equation}\label{4.21}
\Phi_N^{u_1,...,u_k}(t_1,...,t_k)=E\exp\bigg (i\sum_{j=1}^k\big (
(\sum_{l=j}^kt_l)(W_N(u_j)-W_N(u_{j-1}))\big )\bigg ).
\end{equation}
Set
\begin{eqnarray*}
&A_j(N)=\big\{ m:\,[u_{j-1}N]<\te(N)+(m-1)(\te(N)+\tau(N))\\
&<m(\te(N)+\tau(N))\leq [u_jN]\big\},
\end{eqnarray*}
and 
\[
\Psi_N^{u_1,...,u_k}(t_1,...,t_k)=E\exp\bigg (iN^{-1/2}\sum_{j=1}^k
\big ((\sum_{l=j}^kt_l)\sum_{m\in A_j(N)}Y_m\big )\bigg ).
\]
Then similarly to (\ref{4.1}) we show that
\begin{equation}\label{4.22}
|\Phi_N^{u_1,...,u_k}(t_1,...,t_k)-\Psi_N^{u_1,...,u_k}(t_1,...,t_k)|\to 0
\,\,\mbox{as}\,\, N\to\infty.
\end{equation}
Next, similarly to Lemma \ref{lem4.1} we obtain that
\begin{equation}\label{4.23}
|\Psi_N^{u_1,...,u_k}(t_1,...,t_k)-\prod_{j=1}^k\psi_N^{(j)}(t_1,...,t_k)|
\to 0\,\,\mbox{as}\,\, N\to\infty
\end{equation}
where
\[
\psi_N^{(j)}(t_1,...,t_k)=E\exp\big (iN^{-1/2}(\sum_{l=j}^kt_l)
\sum_{m\in A_j(N)}Y_m\big ).
\]
Now in the same way as in (\ref{4.16}) we see that
\begin{equation}\label{4.24}
\psi_N^{(j)}(t_1,...,t_k)\rightarrow\exp\big (-\frac 12\sig^2(u_j-u_{j-1})
(\sum_{l=j}^kt_l)^2\big )\,\,\mbox{as}\,\, N\to\infty
\end{equation}
which together with (\ref{4.22}), (\ref{4.23}) and (\ref{4.17}) yields
(\ref{4.20}).

Next, let $0\leq u_1\leq u\leq u_2\leq 1$ then by Lemma \ref{lem3.7},
\begin{eqnarray}\label{4.25}
& E\big ((W_N(u)-W_N(u_1))^2(W_N(u_2)-W_N(u))^2\big )\\
&\leq\big (E(W_N(u)-W_N(u_1))^4\big )^{1/2}
\big (E(W_N(u_2)-W_N(u))^4\big )^{1/2}\nonumber\\
&\leq\tilde CN^{-2}([uN]-[u_1N])([u_2N]-[uN])\leq\tilde C\big (\frac
{[u_2N]-[u_1N]}{N}\big )^2.\nonumber
\end{eqnarray}
Now, either $u_2-u_1\geq 1/N$ and then the right hand side of (\ref{4.25})
is bounded by $4\tilde C(u_2-u_1)^2$ or $u_2-u_1<1/N$ and then the left hand
side of (\ref{4.25}) is zero. Hence, the left hand side of (\ref{4.25}) is
always bounded by $4\tilde C(u_2-u_1)^2$ and by Ch. 15 of \cite{Bi} the
family $\{\bbP_{W_N},\, N\geq 1\}$ of distributions of $W_N$'s is tight. 
This together with the convergence of finite dimensional distributions of
$W_N$'s established above completes the proof of Theorem \ref{thm2.2}
(cf. Ch. 15 in \cite{Bi}). \qed

\section{Extension to the two linear terms case}\label{sec4+}
\setcounter{equation}{0}

In this section we enhance arguments of Sections \ref{sec3} and \ref{sec4}
in order to derive Theorem \ref{thm2.1+}. 
Set
\begin{eqnarray*}
&R(n)=\prod_{j=0}^\ell X_j(q_j(n))-\prod_{j=0}^\ell
a_j\nonumber\\
&=\sum_{j=0}^\ell a_1\cdots
a_{j-1}(X_j(q_j(n))-a_j)X_{j+1}(q_{j+1}(n))\cdots X_\ell(q_\ell(n)).
\end{eqnarray*}
\begin{lemma}\label{lem4+.1}
There exists $C>0$ such that for all $n\in\bbN$,
\begin{equation}\label{4+.2}
E\big (\sum_{k=0}^nR(k)\big )^2\leq Cn.
\end{equation}
\end{lemma}
\begin{proof}
Relying on (\ref{3.14}) where the summation starts with $j_1,j_2=0$
and estimating $EQ_{j_1j_2}(k_1,k_2)$ essentially by the same argument
as in Lemma \ref{lem3.2} we arrive at (\ref{4+.2}). \qed
\end{proof}

Next, we obtain the 4th moment Gaussian estimate.
\begin{lemma}\label{lem4+.2}
There exists $\tilde C>0$ such that for all $n$ and $m$ satisfying 
$N\geq n>m\geq [N^{1-\gam}]\geq n_0$,
\begin{equation}\label{4+.4}
E\big (\sum_{k=m+1}^nR(k)\big )^4\leq\tilde C(n-m)^2.
\end{equation} 
\end{lemma}
\begin{proof}
Similarly to (\ref{3.34}) and (\ref{3.35}),
\begin{equation}\label{4+.5}
E\big (\sum_{k=m+1}^nR(k)\big)^4\leq\sum_{k_1,k_2,k_3,k_4=m+1}^n
\sum^\ell_{j_1,j_2,j_3,j_4=0}D^{j_1+j_2+j_3+j_4}E|Q_{j_1j_2j_3j_4}
(k_1,k_2,k_3,k_4)|.
\end{equation}
In estimating $|Q_{j_1j_2j_3j_4}(k_1,k_2,k_3,k_4)|$ here we can assume
without loss of generality that $0\leq j_1\leq j_2\leq j_3\leq j_4\leq\ell$.
If $j_1<j_2$ then as in (\ref{3.25}) we still have here that for large
$N$ and $k,l\geq [N^{1-\gam}]$,
\begin{equation}\label{4+.6}
q_{j_2}(l)\geq q_{j_1}(k)+N^\gam -1,
\end{equation}
and so similarly to (\ref{3.26}) taking into account that $k_1,k_2,k_3,k_4>
m\geq [N^{1-\gam}]$ we obtain the estimate (\ref{3.36}) in this
case too. Other estimates of Lemma \ref{lem3.7} hold true here, as well,
since in addition to (\ref{3.3}) and (\ref{4+.6}) we needed there only 
(\ref{3.39}) which is satisfied in the circumstances of Theorem \ref{thm2.1+},
as well. \qed
\end{proof}

Next, we derive a version of Lemma \ref{lem3.5} which holds true under the
conditions of Theorem \ref{thm2.1+}.
\begin{lemma}\label{lem4+.3} There exists $\hat C>0$ such that if $N\geq n
>m\geq [N^{1-\ve}]\geq n_0$ and $n-m\leq\frac 12[N^{1-\ve}]$ for $\ve\in 
(0,\gam)$ then
\begin{equation}\label{4+.7}
|E\big (\sum_{k=m+1}^nR(k)\big )^2-(n-m)(\hat\sig^2-\Xi\prod_{j=2}^\ell
a_j^2)|\leq\hat C
\end{equation}
where $\Xi$ is the same as in (\ref{2.9+++}).
\end{lemma}
\begin{proof} We start with (\ref{3.22}) only now the summation in $j$ should
begin there from 0. For $j\geq 2$ and $k_2>k_1\geq [N^{1-\ve}]$ we still have
 the estimate (\ref{3.24}) while for $j_2>j_1,\, j_2\geq 2$ and $k,l\geq
 [N^{1-\ve}]$ the estimate (\ref{3.26}) holds true though in both cases
 $\ell$ should be replaced by $\ell +1$. Since (\ref{3.27}) remains true 
 also for $j=0$ we obtain (\ref{3.28}) with the summation in $j$ starting
 with 0 and $\ell$ (in the right hand side) replaced by $\ell +1$. Next, 
 (\ref{3.29}) and (\ref{3.30}) remain valid too. Similarly to (\ref{3.29})
 we obtain that for $k_2>k_1\geq [N^{1-\ve}]$,
 \begin{eqnarray}\label{4+.8}
&|EQ_{00}(k_1,k_2)-E\big (X_0(k_2-k_1)-a_0)(X_0(0)-a_0)\big )
E\big(X_1(r(k_2-k_1))X_1(0)\big)\\
&\times\prod_{j=2}^{\ell}a^2_j|\leq 32D^{2(\ell+1)}(\ell+1)
\big ((\ell+1)(\be(\rho_2(N))+\al(\rho_2(N))\big).\nonumber
\end{eqnarray}
Observe that if $\frac 32[N^{1-\ve}]\geq k_1,k_2\geq [N^{1-\ve}]$ then 
$|q_1(k_i)-k_j|\geq\frac 12[N^{1-\ve}],\, i,j=1,2$ and relying on (\ref{3.15})
and (\ref{3.23}) we obtain that
\begin{eqnarray}\label{4+.9}
&|EQ_{01}(k_1,k_2)|\leq 16D^{2(\ell+1)}\bigg ((\ell+1)\be\big(\min(\rho_1(N),
[[N^{1-\ve}]/6])\big)\\
&+\al\big(\min(\rho_1(N),[[N^{1-\ve}]/6])\big)\bigg)\nonumber
\end{eqnarray}
where $\rho_1$ is the same as in (\ref{3.24}). The same estimate holds true
for $|EQ_{10}(k_1,k_2)|$ which together with (\ref{4+.8}), (\ref{4+.9}) and
other estimates mentioned above yield (\ref{4+.7}) similarly to Lemma 
\ref{lem3.5}. \qed
\end{proof}

Next, we enhance arguments of Section \ref{sec4} to make them work in the
situation of Theorem \ref{thm2.1+}. Choose a small $\ve>0$ and a large
$L\geq 4$ so that $L\ve <\gam/4$. Set $\ka(N)=[N^{1-\frac \ve L}],$
$\tau(N)=[N^{1-\ve}],$ $\te(N)=[N^{1-L\ve}]$ and $\mu(N)=[\frac {(r-1)\ka(N)
+p}{\tau(N)+\te(N)}]$ recalling that $q_1(n)=rn+p$ and assuming that 
$N\geq\exp(2/\ve)$ which ensures that $\mu(N)\geq 1$. Using the notation
$q_1^{(l)}(n)=q_1(q_1^{(l-1)}(n)),$ $q_1^{(1)}=q_1$ for iterates of $q_1$
define
\begin{eqnarray*}
&L_{kl}(N)=q_1^{(l)}\big(\ka(N)+(k-1)(\tau(N)+\te(N))\big)\,\,\mbox{and}\\
&\tilde L_{kl}(N)=q_1^{(l)}\big(\ka(N)+\tau(N)+(k-1)(\tau(N)+\te(N))\big).
\end{eqnarray*}
Introduce the sets of integers
\[
\Gam_{kl}(N)=\{ n:\, L_{kl}(N)\leq n<\tilde L_{kl}(N)\}\,\,\mbox{and}\,\,
\tilde\Gam_{kl}(N)=\{ n:\,\tilde L_{kl}(N)\leq n<L_{k+1,l}(N)\}
\]
where $k=1,2,...,\mu(N)$ and $l=1,2,...,\nu_k(N)$ with $\nu_k(N)=\max\{
l:\, L_{k+1,l}(N)\leq N\}$. The block sequences 
$\{\Gam_{kl}(N)\}_{l=1}^{\nu_k(N)}$ and 
$\{\tilde\Gam_{kl}(N)\}_{l=1}^{\nu_k(N)}$ will play the same role as the
blocks $\Gam_k(N)$ and $\tilde\Gam_k(N)$ in Section \ref{sec4}.

Set
\begin{eqnarray*}
&Y_{kl}=\sum_{n\in\Gam_{kl}(N)}R(n),\,\, Z_{kl}=\sum_{n\in\tilde\Gam_{kl}(N)}
R(n),\,\, Y_k=\sum_{1\leq l\leq\nu_k(N)}Y_{kl}\\
&\Phi_N(t)=E\exp\big(\frac {it}{\sqrt N}\sum_{n=0}^NR(n)\big)\,\,\mbox{and}\,\,
\Psi_N(t)=E\exp\big(\frac {it}{\sqrt N}\sum_{1\leq k\leq\mu(N)}Y_k\big).
\end{eqnarray*}
Similarly to (\ref{4.1}) we obtain from (\ref{4+.2}) and (\ref{4+.4})
that for some constant $\check C>0$ and all $t$ and $N$,
\begin{equation}\label{4+.10}
|\Phi_N(t)-\Psi_N(t)|\leq\check C|t|\big(N^{-\frac {\ve}{2L}}+
N^{(1-\frac L2)\ve}\ln N+N^{-\frac {(L-1)\ve}{2L}}\ln N\big).
\end{equation}

Set 
\[
\psi_N^{(k)}(t)=E\exp\big(\frac {it}{\sqrt N}Y_k\big).
\]
\begin{lemma}\label{lem4+.4} For any $t$ and each small $\ve>0$ there
exists $K_\ve(t)>0$ such that for all $N\geq\exp(2/\ve)$,
\begin{equation}\label{4+.11}
|\Psi_N(t)-\prod_{1\leq k\leq\mu(N)}\psi_N^{(k)}(t)|\leq K_\ve(t)N^{-\frac 
\ve 2\sqrt N}.
\end{equation}
\end{lemma}
\begin{proof}
The argument goes on, essentially, in the same way as in Lemma \ref{lem4.1}.
Namely, we set
\[
\hat Y_{kl}=Y_{kl}+\tau(N)\prod_{j=0}^\ell a_j,\,\,\hat Y_k=\sum_{1\leq l\leq
\nu_k(N)}\hat Y_{kl}
\]
and proceed as in (\ref{4.3})--(\ref{4.9}). Next, we write
\begin{equation}\label{4+.12}
\big (\sum_{1\leq m\leq\nu_k(N)}\hat Y_{km}\big )^{l_k}=\sum_{\sig^{(k)}}
\be^{(k)}_{\sig^{(k)}}I^{(k)}_{\sig^{(k)}}J^{(k)}_{\sig^{(k)}}
\end{equation}
where
\[
I^{(k)}_{\sig^{(k)}}=\prod_{1\leq m\leq\nu_k(N)}\prod_{n\in\Gam_{km}(N)}
X_0^{\sig^{(k)}_n}(n)X_1^{\sig^{(k)}_n}(q_1(n)),
\]
\[
J^{(k)}_{\sig^{(k)}}=\prod_{1\leq m\leq\nu_k(N)}\prod_{n\in\Gam_{km}(N)}
\prod_{j=2}^\ell X_j^{\sig^{(k)}_n}(q_j(n)),
\]
$\be^{(k)}_{\sig^{(k)}}$ are $l_k$-nomial coefficients and $\sig^{(k)}=
\big (\sig^{(k)}_n,\, n\in\Gam_{km}(N),\, m=1,...,\nu_k(N)\big )$ satisfies
\begin{equation}\label{4+.13}
\sig^{(k)}_n\geq 0\,\,\mbox{and}\,\,
\sum_{1\leq m\leq\nu_k(N)}\sum_{n\in\Gam_{km}(N)}\sig^{(k)}_n=l_k\leq n(N)
\end{equation}
with $n(N)$ defined by (\ref{4.5}). Then we obtain the estimate (\ref{4.12})
with $\mu(N)$ in place of $m(N)$ and
\begin{equation}\label{4+.14}
H_{l_1,...,l_{\mu(N)}}(t,N)=\big\vert E\prod_{k=1}^{\mu(N)}I^{(k)}_{\sig^{(k)}}
J^{(k)}_{\sig^{(k)}}-\prod_{k=1}^{\mu(N)}EI^{(k)}_{\sig^{(k)}}
J^{(k)}_{\sig^{(k)}}\big\vert.
\end{equation}
Since all $n\in\Gam_{km}(N)$ satisfy $n\geq N^{1-\frac \ve L}$ we obtain
from (\ref{3.3}) and (\ref{3.23}) similarly to (\ref{4.13}) that
\begin{eqnarray}\label{4+.15}
&\big\vert E\prod_{k=1}^{\mu(N)}I^{(k)}_{\sig^{(k)}}J^{(k)}_{\sig^{(k)}}-
E(\prod_{k=1}^{\mu(N)}I^{(k)}_{\sig^{(k)}})E(\prod_{k=1}^{\mu(N)}
J^{(k)}_{\sig^{(k)}})\big\vert\\
&\leq 2D^{(\ell+1)n(N)\mu(N)}\big ((\ell+1)n(N)\mu(N)\be(\rho_7(N))+
2\al(\rho_7(N))\big ),\nonumber
\end{eqnarray}
where 
\[
\rho_7(N)=[\frac 13([N^{1-\frac \ve L}]-[N^{(1-\gam)(1-\frac \ve L)}])],
\]
and
\begin{eqnarray}\label{4+.16}
&\big\vert EI^{(k)}_{\sig^{(k)}}J^{(k)}_{\sig^{(k)}}-
EI^{(k)}_{\sig^{(k)}})EJ^{(k)}_{\sig^{(k)}})\big\vert\\
&\leq 2D^{(\ell+1)n(N)}\big ((\ell+1)n(N)\be(\rho_7(N))+
2\al(\rho_7(N))\big ).\nonumber
\end{eqnarray}

Observe that if $n_1\in\Gam_{k_1m_1}(N)$ and $n_2\in\Gam_{k_2m_2}(N)$
with either $k_1\ne k_2$ or $m_1\ne m_2$ then $|n_1-n_2|>\te(N)=[N^{1-L\ve}]$.
Thus using (\ref{3.25}), (\ref{3.40}) and applying (\ref{3.3}) no more than 
$2\sum_{1\leq k\leq\mu(N)}\nu_k(N)$ times we obtain that
\begin{eqnarray}\label{4+.17}
&\big\vert E\prod_{k=1}^{\mu(N)}J^{(k)}_{\sig^{(k)}}-
\prod_{k=1}^{\mu(N)}EJ^{(k)}_{\sig^{(k)}}\big\vert
\leq 4(\sum_{1\leq k\leq\mu{N}}\nu_k(N))D^{(\ell-1)n(N)\mu(N)}\\
&\times \big ((\ell-1)n(N)\mu(N)\be(\rho_8(N))+
2\al(\rho_8(N))\big )\nonumber
\end{eqnarray}
where 
\[
\rho_8(N)=[\frac 13\min(\te(N),\rho_7(N))].
\]
By our construction if $n\in\Gam_{km}(N)$ then $q_1(n)\in\Gam_{k,m+1}(N)$,
and so we can represent $I^{(k)}_{\sig^{(k)}}$ in the form
\[
I^{(k)}_{\sig^{(k)}}=\prod_{1\leq m\leq\nu_k(N)}\,\prod_{n\in\Gam_{km}(N)}
X_0^{\eta^{(k)}_n}(n)X_1^{\zeta^{(k)}_n}(n)
\]
which together with (\ref{3.3}) and the above argument that $|n_1-n_2|>\te(N)$
for $n_j\in\Gam_{k_jm_j},\, j=1,2$ from different blocks yields that
\begin{eqnarray}\label{4+.18}
&\big\vert E\prod_{k=1}^{\mu(N)}I^{(k)}_{\sig^{(k)}}-\prod_{k=1}^{\mu(N)}
\prod_{1\leq m\leq\nu_k(N)}E\prod_{n\in\Gam_{km}(N)}
X_0^{\eta^{(k)}_n}(n)X_1^{\zeta^{(k)}_n}(n)\big\vert\\
&\leq 8(\sum_{1\leq k\leq\mu{N}}\nu_k(N))D^{2n(N)\mu(N)}
\big (n(N)\mu(N)\be([\frac 13\te(N)])+\al([\frac 13\te(N)])\big ).\nonumber
\end{eqnarray}
Similarly,
\begin{eqnarray}\label{4+.19}
&\big\vert EI^{(k)}_{\sig^{(k)}}-
\prod_{1\leq m\leq\nu_k(N)}E\prod_{n\in\Gam_{km}(N)}
X_0^{\eta^{(k)}_n}(n)X_1^{\zeta^{(k)}_n}(n)\big\vert\\
&\leq 8\nu_k(N)D^{2n(N)}
\big (n(N)\be([\frac 13\te(N)])+\al([\frac 13\te(N)])\big ).\nonumber
\end{eqnarray}
Collecting (\ref{4+.14})--(\ref{4+.19}) we obtain that
\begin{eqnarray}\label{4+.20}
&H_{l_1,...,l_{\mu(N)}}(t,N)\leq 28\big (\sum_{1\leq k\leq\mu{N}}\nu_k(N))
D^{(\ell +1)n(N)\mu(N)}\mu(N)\\
&\times\big ((\ell+1)n(N)\be(\rho_8(N))+
2\al(\rho_8(N))\big ).\nonumber
\end{eqnarray}

Observe that for each $k$,
\begin{equation}\label{4+.21}
\sum_{\sig^{(k)}}\be^{(k)}_{\sig^{(k)}}\leq (\tau(N)\sum_{l=0}^{\nu_k(N)}r^l)^
{l_k}\leq (N^{1-\ve}r^{\nu_k(N)})^{n(N)}\leq N^{(1-(1-\frac 1L)\ve)
N^{\frac 12+\ve}}
\end{equation}
since, clearly,
\begin{equation}\label{4+.22}
\nu_k(N)\leq\frac {\ve\ln N}{L\ln r}.
\end{equation}
In addition, we see by (\ref{4+.22}) that
\begin{equation}\label{4+.23}
\sum_{1\leq l_1,...,l_{\mu (N)}\leq n(N)}\prod_{k=1}^{\mu(N)}\frac
{|N^{\frac 12-\ve}tr^{\nu_k(N)}|^{l_k}}{l_k!}\leq\exp(N^{\frac 12-
\ve(1-\frac 1L)}|t|\mu(N)).
\end{equation}
Employing (\ref{4.3})--(\ref{4.9}) and (\ref{4.12}) with $\mu(N)$ in place 
of $m(N)$ and $H_{l_1,...,l_{\mu(N)}}$ given by (\ref{4+.14}) together with
(\ref{4+.20})--(\ref{4+.23}) we arrive at (\ref{4+.11}).   \qed
\end{proof}

Next, in order to complete the proof of Theorem \ref{thm2.1+} in the same way
as at the end of Section \ref{sec4} proceeding via (\ref{4.16})--(\ref{4.19})
we observe that by (\ref{4+.4}) if $M_{km}(N)$ denotes the number of
integers in $\Gam_{km}(N)$ then 
\begin{eqnarray}\label{4.+24}
&E\big (\sum_{m=1}^{\nu_k(N)}\sum_{n\in\Gam_{km}(N)}R(n)\big )^4\leq
(\nu_k(N))^3\sum_{m=1}^{\nu_k(N)}E(\sum_{n\in\Gam_{km}(N)}R(n))^4\\
&\leq\tilde C(\nu_k(N))^3\sum_{m=1}^{\nu_k(N)}(M_{km}(N))^2\leq\tilde C
(\nu_k(N))^3\big (\sum_{m=1}^{\nu_k(N)}M_{km}(N)\big )^2\nonumber
\end{eqnarray}
which in view of (\ref{4+.22}) is still sufficient for the estimate of
the form (\ref{4.16}).

Finally, we show that 
\begin{equation}\label{4+.25}
A_k(N)=\big\vert E\big (\sum_{m=1}^{\nu_k(N)}\sum_{n\in\Gam_{km}(N)}R(n)
\big )^2-\hat\sig^2M_k(N)\big\vert\leq C\nu_k(N)
\end{equation}
where $M_k(N)=\sum_{1\leq m\leq\nu_k(N)}M_{km}(N)$ and $C>0$ does not depend
on $N$ and $k$. Indeed, by (\ref{4+.7}),
\begin{equation}\label{4+.26}
\big\vert\sum_{m=1}^{\nu_k(N)}E(\sum_{n\in\Gam_{km}(N)}R(n))^2-M_k(N)
(\hat\sig^2-\Xi\prod_{j=2}^\ell a_j^2)\big\vert\leq\hat C\nu_k(N).
\end{equation}
Observe that if $m_2-m_1\geq 2$ then for any $n_1\in\Gam_{km_1}$ and 
$n_2\in\Gam_{km_2}$ we have that $n_2-q_1(n_1)\geq [N^{1-L\ve}].$ Thus using
(\ref{3.15}), (\ref{3.23}) and (\ref{3.14}) (the latter with the summation
 starting with $j_1,j_2=0$) we obtain for such $n_1$ and 
$n_2$ that
\begin{equation}\label{4+.27}
|R(n_1)R(n_2)|\leq 16D^{4(\ell+1)}(\ell+1)^2\big ((\ell+1)\be(\rho_9(N))+
\al(\rho_9(N))\big )
\end{equation}
where
\[
\rho_9(N)=[\frac 13\min([N^{1-L\ve}],\, N^{1-\frac \ve L}-[N^{(1-\frac \ve L)
(1-\gam)}])].
\]

It remains to estimate
\begin{equation}\label{4+.28}
B_{km}=E\big ((\sum_{n\in\Gam_{km}(N)}R(n))(\sum_{\tilde n\in\Gam_{k,m+1}(N)}
R(\tilde n))\big )=\sum_{n\in\Gam_{km}(N),\tilde n\in\Gam_{k,m+1}(N)}
ER(n)R(\tilde n).
\end{equation}
Using (\ref{3.15}), (\ref{3.22}) and (\ref{3.23}) we obtain that 
\begin{eqnarray}\label{4+.29}
&|B_{km}-a_0\sum_{n\in\Gam_{km}(N),\tilde n\in\Gam_{k,m+1}(N)}EQ_{01}(\tilde n,
n)|\\
&\leq 48D^{4(\ell+1)}(\ell+1)^2\big ((\ell+1)\be(\rho_9(N))+\al(\rho_9(N))
\big ).\nonumber
\end{eqnarray}
In the same way as in (\ref{4+.8}) it follows that for $n\in\Gam_{km}(N)$ and
 $\tilde n\in\Gam_{k,m+1}(N)$,
 \begin{eqnarray}\label{4+.30}
 &|EQ_{01}(\tilde n,n)-a_1E(X_0(\tilde n)-a_0)(X_1(q_1(n)-a_1)
 \prod_{\ell=2}^\ell a_j^2|\\
 &\leq 32D^{2(\ell+1)}(\ell+1)\big ((\ell+1)\be(\rho_2(N))+\al(\rho_2(N))\big ).
 \nonumber\end{eqnarray}
 Furthermore,
  \begin{eqnarray}\label{4+.31}
  &\sum_{n\in\Gam_{km}(N),\tilde n\in\Gam_{k,m+1}(N)}E(X_0(\tilde n)-a_0)
  (X_1(q_1(n))-a_1)\\
  &= \sum_{n\in\Gam_{km}(N)}\big (\sum_{\tilde n\geq q_1(n)}
  E(X_0(\tilde n-q_1(n))-a_0)(X_1(0)-a_1)\nonumber\\
  &+\sum_{\tilde n< q_1(n)}E(X_0(0)-a_0)(X_1(q_1(n)-\tilde n)-a_1)-V_1-V_2
  \nonumber\end{eqnarray}
  where by (\ref{2.1}) and (\ref{3.3}),
  \begin{eqnarray}\label{4+.32}
  &|V_1|+|V_2|\leq\sum_{n\in\Gam_{km}(N)}\big (\sum_{\tilde n\geq
  \tilde L_{k,m+1}(N)}\exp(c(\tilde n-q_1(n)))\\
  &+\sum_{\tilde n<L_{k,m+1}(N)}\exp(c(q_1(n)-\tilde n))\big )\leq\tilde 
  {\tilde C}\nonumber
  \end{eqnarray}
  for some $c,\tilde {\tilde C}>0$ independent of $N,k$ and $m$. Collecting
  (\ref{4+.26})--(\ref{4+.32}) we obtain (\ref{4+.25}). In view of 
  (\ref{4+.22}) this enables us to complete the proof of Theorem \ref{thm2.1+}
  in the same way as at the end of Section \ref{sec4}.  \qed

\section{Concluding remarks}\label{sec5}\setcounter{equation}{0}

The condition (\ref{2.4}) was crucial for our proof of Theorem \ref{thm2.1}
since its, essentially, equivalent form (\ref{3.25}) arranges $q_j(n),\,
j=1,...,\ell$ for big $n$ into $\ell$ sets separated by large gaps
which was necessary in our splitting arguments. This property is
lost when more than one of $q_j$'s are linear but, still, the block 
sequences construction of Section \ref{sec4+} enabled us to carry out the
proof of Theorem \ref{thm2.1+} for two linear terms. 
Lemmas \ref{lem4+.1}--\ref{lem4+.3} still can be carried out when more than
two $q_j$'s are linear but it is not clear how to make an appropriate block 
sequences construction in this case, for instance, when $q_1(n)=n,\, q_2(n)=2n,
\, q_3(n)=3n$ and $\ell=3$. 
Probably, in a special algebraic situation, for instance, when
$X_j(n)=X(n)=f(T^nx)$ with $T$ being a hyperbolic automorphism or an
expanding (algebraic) endomorphism of a torus, the Fourier analysis
technique in the spirit of \cite{Fu} may still lead to a SPLIT in
the form of Theorems \ref{thm2.1}--\ref{thm2.1+}. Nevertheless, for more
general stationary processes $X(n)$ it is not clear whether a Theorems
\ref{thm2.1}--\ref{thm2.1+} type result holds true for expressions of the form
\begin{equation}\label{5.1}
N^{-1/2}\sum_{0\leq n\leq N}\big (X(n)X(2n)X(3n)-(EX(0))^3\big ).
\end{equation}
On the other hand, if $X(0),X(1),X(2),...$ are i.i.d. random
variables such results can be easily proved. Namely, let
$q_1=1<q_2<...<q_\ell$ be some prime numbers and set $EX^2(0)=b^2$
assuming for simplicity that $EX(0)=0$. Then as $N\to\infty$,
\begin{equation}\label{5.2}
W_N=N^{-1/2}\sum_{0\leq n\leq N}X(q_1n)X(q_2n)\cdots X(q_\ell n)
\end{equation}
converges in distribution to the centered normal random variable
with the variance $\sigma^2=b^{2\ell}$. Indeed, let $1\leq
k_1<k_2<...<k_{m_N}\leq N$ be all integers which are not divisible
by any of $q_j$'s, $j\geq 2$. Then we can define disjoint sets
$A_{k_l},\, l=1,...,m_N$ so that $A_{k_l}\subset\{1,...,N\}$ and any
$n\in A_{k_l}$ is obtained from $k_l$ by multiplication by some of
$q_j$'s. It is clear that the number $r_\ell(N)$ of elements of each
$A_{k_l}$ does not exceed $\log_2N$. Set
\begin{equation}\label{5.3}
S_N(l)=\sum_{n\in A_{k_l}}X(q_1n)X(q_2n)\cdots X(q_\ell n).
\end{equation}
Then $S_N(l),\, l=1,2,...,m_N$ are independent random variables with
zero mean and the variance $r_l(N)b^2$. Applying the standard
central limit theorem for triangular arrays (see, for instance,
\cite{Sh}) to 
\begin{equation}\label{5.4}
W_N=N^{-1/2}\sum_{0\leq l\leq m_N}S_N(l)
\end{equation}
and taking into account that $\sum_{0\leq l\leq m_N}r_l(N)=N$ we obtain the
required result. If $EX(0)\ne 0$ then this method still works using
the representation (\ref{3.10}) for computation of variances.

Observe that, in principle, we could ask whether under appropriate
conditions our results could be extended to continuous time processes trying
to obtain central limit theorems for integrals 
\[
\int_0^TX_1(q_1(t))X_2(q_2(t))\cdots X_\ell(q_\ell(t))dt
\]
 in place of sums. Nevertheless, the answer to this
question is not clear yet and the approach of this paper does not seem to work
in this case.

Another result which can be derived for i.i.d. bounded random variables\newline
$X(0), X(1), X(2),...$ is a corresponding sum-product large deviations 
(SPLAD) theorem. Namely, we are interested in the asymptotic behavior of
\begin{equation}\label{5.5}
Q_N(U)=\frac 1N\log P\{\frac 1NS_N\in U\}
\end{equation}
as $N\to\infty$ where $S_N=\sum^N_{n=0}X(q_1(n))\cdots X(q_\ell(n))$ and
$U\subset\bbR$. Here, $q_1(n),...,q_\ell(n)$ are nonnegative strictly 
increasing functions taking on integer values on the integers and such 
that for some $\gam\in(0,1)$ and $n_0\in\bbZ$ we have 
\begin{equation}\label{5.6}
q_{j+1}([n^\gam])>q_j(n)\quad\mbox{for all}\,\,\, n\geq n_0.
\end{equation}
Let $M_N(t)=E\exp(tS_N)$ be the moment generating function of $S_N$. It is
well known (see, for instance, Theorem 2.3.6 in \cite{DZ}) that if the limit
\begin{equation}\label{5.7}
\eta(t)=\lim_{N\to\infty}\frac 1N\log M_N(t)
\end{equation}
exists and it is differentiable in $t$ then
\begin{equation}\label{5.8}
\limsup_{N\to\infty}Q_N(F)\leq -\inf_{x\in F}\La(x)
\end{equation}
for any closed set $F\subset\bbR$ and 
\begin{equation}\label{5.9}
\liminf_{N\to\infty}Q_N(G)\geq -\inf_{x\in G}\La(x)
\end{equation}
for any open set $G\subset\bbR$ where
\[
\La(x)=\sup_t(tx-\eta(t))
\]
is the Legendre transform of $\eta$.

Set $\tilde S_N=\sum_{N\geq n\geq [N^\gam]}X(q_1(n))\cdots X(q_\ell(n))$
and $\tilde M_N(t)=E\exp(t\tilde S_N)$. Then 
\begin{equation}\label{5.10}
\tilde M_N(t)\exp(-D^\ell N^\gam)\leq M_N(t)\leq\tilde M_N(t)
\exp(D^\ell N^\gam),
\end{equation}
where we assume that $|X(0)|\leq D$ a.s., and so
\begin{equation}\label{5.11}
\lim_{N\to\infty}\frac 1NM_N(t)=\lim_{N\to\infty}\frac 1N\log\tilde M_N(t)
\end{equation}
whenever  one of these limits exists. Set $m_n(t)=E\exp\big(t\prod_{j=1}^\ell
X(q_j(n))\big)$. 
By (\ref{5.6}) and the strict monotonicity of the functions $q_j(n)$ it
follows that the terms\newline
 $\exp\big(t\prod_{j=1}^\ell X(q_j(n))\big)$ are
independent for different $n$, and so
\begin{equation}\label{5.12}
\tilde M_N(t)=\prod_{n=[N^\gam]}^Nm_n(t).
\end{equation}
Next, by (\ref{5.6}) the factors in the product appearing in the definition of
 $m_n(t)$ with $n\geq [n_0^\gam]$ are independent, and so for such $n$,
\begin{eqnarray}\label{5.13}
&m_n(t)=E\sum_{k=0}^\infty\frac {t^k}{k!}\prod_{j=1}^\ell X^k(q_j(n))\\
&=\sum_{k=0}^\infty\frac {t^k}{k!}\prod_{j=1}^\ell EX^k(q_j(n))=
\sum_{k=0}^\infty\frac {t^k}{k!}\prod_{j=1}^\ell EX^k(0)=m_{[n_0^\gam]}(t).
\nonumber\end{eqnarray}
Thus, we obtain
\begin{equation}\label{5.14}
\eta(t)=\lim_{N\to\infty}\frac 1N\log\tilde M(t)=\log m_{[n_0^\gam]}(t),
\end{equation}
which is, clearly, differentiable
in $t$ since $X(k)$'s are bounded, and so (\ref{5.8}) and (\ref{5.9}) follow.
SPLAD in other situations will be treated in another paper.
 
For i.i.d. $X(j),\, j=0,1,2,...$ it is easy to prove the existence of a
differentiable limit $\eta(t)$ in (\ref{5.7}) also for moment generating
functions $M_N(t)$ of the sums
\[
S_N=\sum_{0\leq n\leq N}X(q_1n)X(q_2n)\cdots X(q_\ell n),
\]
where $q_i,\, i=1,...,\ell$ are primes as in (\ref{5.2}), by using the sets 
$A_{k_l}$ and partial sums $S_N(l)$ appearing in (\ref{5.3}). 

In conclusion, remark that using the thermodynamic formalism and decay of 
correlations results for random transformations from \cite{Ki1} and \cite{Ki2}
we can obtain the corresponding (quenched or fiberwise) SPLIT for random 
subshifts of finite type, random expanding transformations and for Markov 
chains with random transitions.

\end{document}